\documentclass[leqno,final]{siamltex}
\usepackage{graphicx} 
\usepackage{amsmath,amstext,amssymb,bm}
\usepackage{leftidx}

\usepackage{xcolor} 
\usepackage{soul} 
\usepackage{tikz}
\usetikzlibrary{shapes,arrows}
\usepackage{mathrsfs}
\usepackage{epstopdf}
\usepackage{color}
\usepackage{multirow}
\usepackage{tabularx}
\usepackage[shortlabels]{enumitem}

\numberwithin{equation}{section}
\newtheorem{remark}{Remark}[section]

\allowdisplaybreaks[4]

\renewcommand{\div}{\mbox{\rm div\,}}

\newcommand{\cK}{\mathcal{K}}

\newcommand{\cF}{\mathcal{F}}

\newcommand{\mP}{\mathbb{P}}
\newcommand{\mH}{\mathbb{H}}
\newcommand{\mE}{\mathbb{E}}
\newcommand{\mV}{\mathbb{V}}

\newcommand{\eps}{\epsilon}

\newcommand{\Ome}{\Omega}
\newcommand{\p}{\partial}
\newcommand{\nab}{\nabla}
\newcommand{\vu}{{\bf u}}

\newcommand{\vP}{{\bf P}}
\newcommand{\vf}{{\bf f}}
\newcommand{\vH}{{\bf H}}

\newcommand{\vL}{{\bf L}}
\newcommand{\vv}{{\bf v}}
\newcommand{\ve}{{\bf e}}
\newcommand{\vG}{{\bf G}}

\newcommand{\vQ}{\mathbf{Q}}

\newcommand{\G}{\pmb{\mathcal{G}}}
\newcommand{\pphi}{\pmb{\phi}}

\newcommand{\vA}{{\bf A}}

\begin{document}
	
	\title{Higher order time discretization method for the stochastic Stokes equations with multiplicative noise}
	\markboth{LIET VO}{THE STOCHASTIC STOKES EQUATIONS}

	\author{Liet Vo\thanks{Department of Mathematics, Statistics and Computer Science, The University of Illinois at Chicago, Chicago, IL 60607, U.S.A. (lietvo@uic.edu).}}

	\maketitle
	
	\begin{abstract} In this paper, we propose a new approach for the time-discretization of the incompressible stochastic Stokes equations with multiplicative noise. Our new strategy is based on the classical Milstein method from stochastic differential equations. We use the energy method for its error analysis and show a strong convergence order of at most $1$ for both velocity and pressure approximations. The proof is based on a new H\"older continuity estimate of the velocity solution. While the errors of the velocity approximation are estimated in the standard $L^2$- and $H^1$-norms, the pressure errors are carefully analyzed in a special norm because of the low regularity of the pressure solution. In addition, a new interpretation of the pressure solution, which is very useful in computation, is also introduced. Numerical experiments are also provided to validate the error estimates and their sharpness.
	\end{abstract}
	
	
	
	

	\section{Introduction}\label{sec-1}
	In this paper, we consider the following incompressible time-dependent stochastic Stokes equations with multiplicative noise:
	\begin{subequations}\label{eq1.1}
		\begin{alignat}{2} \label{eq1.1a}
			d\vu &=\bigl[\nu\Delta \vu-\nabla p + \vf\bigr]\,dt + {\vG}(\vu) \, dW(t)  &&\qquad\mbox{a.s. in}\, D_T:=(0,T)\times D,\\
			\div \vu &=0 &&\qquad\mbox{a.s. in}\, D_T,\label{eq1.1b}\\
			\vu(0)&= \vu_0 &&\qquad\mbox{a.s. in}\, D,\label{eq1.1d}
		\end{alignat}
	\end{subequations}
	where $D = (0, L)^d\subset \mathbb{R}^d \, (d=2,3)$ represents a period of the periodic domain in $\mathbb{R}^d$, $\vu$ and $p$ stand for respectively the velocity field and the pressure of the fluid,
	$\vG$ is an operator-valued random field, $\{W(t); t\geq 0\}$ denotes an ${{\mathbb{R}}}$-valued
	Wiener process and $\vf$ is a body force function (see Section \ref{sec2} for their precise definitions).  Here we seek periodic-in-space solutions $(\vu,p)$ with period $L$, that is, 
	$\vu(t,{\bf x} + L{\bf e}_i) = \vu(t,{\bf x})$ and $p(t,{\bf x}+L{\bf e}_i)=p(t,{\bf x})$ 
	almost surely 
	and for any $(t, {\bf x})\in (0,T)\times \mathbb{R}^d$  and $1\leq i\leq d$, where 
	$\{\bf e_i\}_{i=1}^d$ denotes the canonical basis of $\mathbb{R}^d$.
	
	The stochastic Stokes system is one of the most well-known equations in fluid mechanics due to its broad applications in modeling micro fluids, where microscopic fluctuations are pertinent contributions to the flow of fluid dynamics and inertial issues are generally careless \cite{CHP2012, squires2005, dellacherie1385}. Unlike the deterministic counterpart, the system \eqref{eq1.1} is a nonlinear SPDE due to the nonlinearity of the diffusion coefficient $\vG$.	On the PDE side of \eqref{eq1.1}, we refer the reader to intensive works in establishing the existence and uniqueness of the solutions in \cite{LRS03, bensoussan1973equations,flandoli1995martingale,menaldi2002stochastic}. The numerical analysis of the stochastic Stoke equations in \eqref{eq1.1} has been gradually developed in the last decade. We can list a few works in the literature until now. In \cite{CHP2012}, the authors considered splitting strategies for the time discretization of \eqref{eq1.1} in the case of divergence-free noise. Later, the authors in \cite{FL2020} extended the splitting strategies for the same equations but they considered non-divergence-free noise. In \cite{FQ, FPL2021}, the authors studied an implicit Euler-Maruyama method for the time-discretization and mixed finite element methods for the spatial discretization of \eqref{eq1.1}. On the stochastic Navier-Stokes equations, the convective nonlinear term of the Navier-Stokes makes the error analysis completely different from its deterministic counterpart and of course the stochastic Stokes equations. In \cite{CHP2012}, the authors employed the standard Euler-Maruyama method for time-discretization of the stochastic Navier-Stokes equations with multiplicative noise. They addressed the difficulty, which was caused by the nonlinear term, in obtaining a strong convergence for the semi-discrete solutions. 
	In addition, due to their low order H\"older continuity estimate, they just obtained a convergence order of $\frac14$ for the semi-discrete velocity approximation. Later, the authors in \cite{Breit} improved the H\"older continuity estimate for the velocity solution and they established a convergence order of $\frac12$ for the semi-discrete velocity approximation.

	To the best of our knowledge, all the time discretization methods for the stochastic Stokes and Navier-Stokes equations in the literature until now are just providing a convergence order of at most $\frac12$. Therefore, the first goal of this paper is to propose a new approach that is based on the classical Milstein method \cite{mil1975approximate} in the stochastic differential equations (SDEs) for the time discretization of \eqref{eq1.1}. Our new approach (see Algorithm 1) can produce a convergence order of at most $1$ for both velocity and pressure approximations. In addition, we employ the energy method for the error analysis of the semi-discrete velocity solution. The proof is based on a new (higher-order) H\"older continuity estimate of the velocity solution (see Lemma \ref{lemma2.4}). 
	
	Unlike the velocity, the pressure solution is a distribution only \cite{LRS03}, and therefore it is unclear how to obtain an error estimate for the pressure solution $p$  of \eqref{eq1.1}. Recently, the papers \cite{FQ, FPL2021, FL2020} have come with a new interpretation of $p$ as the distributional time derivative of some time-averaged pressure $P$ (see Theorem \ref{thm 2.2}). The authors proved that the semi-discrete pressure approximation converged to the time-averaged pressure $P$ with an order of $\frac12$. We notice that the convergence of the pressure approximation is weak. Therefore, another goal of this paper is to introduce a new interpretation of the pressure solution and then a stronger error estimate of the pressure approximation. 
	
	The remainder of this paper is organized as follows. In Section \ref{sec2}, we introduce the mathematical notations that will be used in this paper. We also state assumptions of the diffusion functional $\vG$. In Subsection \ref{sec2.1}, we recall the variational solution of \eqref{eq1.1} and prove the higher order H\"older continuity estimate for the velocity solution in Lemma \ref{lemma2.4}. In Subsection \ref{sec2.2}, we present a new approach to the pressure solution $p$ by defining a new ``simulate" pressure solution. In Section \ref{section_semi}, we state the proposed method in Algorithm 1 and prove the important technical Lemma \ref{lemma3.1}. Next, in Subsection \ref{sec3.2}, we establish the convergence order $O(k^{1-\eps})$ for the velocity and pressure approximations. In Section \ref{section_fullydiscrete}, we present the standard mixed finite element approximation of Algorithm 1 and establish the global error estimates for velocity and pressure approximations. Finally, in Section \ref{sec_compute}, we present a few numerical tests to verify our theoretical results.

	\section{Preliminaries}\label{sec2}
	Standard function and space notation will be adopted in this paper. 
	Let $\vH^1_0(D)$ denote the subspace of $\vH^1(D)$ whose ${\mathbb R}^d$-valued functions have zero trace on $\p D$, and $(\cdot,\cdot):=(\cdot,\cdot)_D$ denote the standard $L^2$-inner product, with induced norm $\Vert \cdot \Vert$. We also denote ${\bf L}^p_{per}(D)$ and ${\bf H}^{k}_{per}(D)$ as the Lebesgue and Sobolev spaces of the functions that are periodic and have vanishing mean, respectively. 
	Let $(\Omega,\cF, \{\cF_t\},\mP)$ be a filtered probability space with the probability measure $\mP$, the 
	$\sigma$-algebra $\cF$ and the continuous  filtration $\{\cF_t\} \subset \cF$. For a random variable $v$ 
	defined on $(\Omega,\cF, \{\cF_t\},\mP)$,
	${\mathbb E}[v]$ denotes the expected value of $v$. 
	For a vector space $X$ with norm $\|\cdot\|_{X}$,  and $1 \leq p < \infty$, we define the Bochner space
	$\bigl(L^p(\Omega,X); \|v\|_{L^p(\Omega,X)} \bigr)$, where
	$\|v\|_{L^p(\Omega,X)}:=\bigl({\mathbb E} [ \Vert v \Vert_X^p]\bigr)^{\frac1p}$.
	We also define 
	\begin{align*}
		{\mathbb H} &:= \bigl\{{\bf v}\in  \vL^2_{per}(D) ;\,\div {\bf v}=0 \mbox{ weakly in }D\, \bigr\}\, , \\
		{\mathbb V} &:=\bigl\{{\bf v}\in  \vH^1_{per}(D) ;\,\div {\bf v}=0 \mbox{ weakly in }D \bigr\}\, .
	\end{align*}
	
	We recall from \cite{Girault_Raviart86} that the (orthogonal) Helmholtz projection 
	${\bf P}_{{\mathbb H}}: \vL^2_{per}(D) \rightarrow {\mathbb H}$ is defined 
	by ${\bf P}_{{\mathbb H}} {\bf v} = \pmb{\eta}$ for every ${\bf v} \in \vL^2_{per}(D)$, 
	where $(\pmb{\eta}, \xi) \in {\mathbb H} \times H^1_{per}(D)/\mathbb{R}$ 
	We define the Stokes operator ${\bf A} := -{\bf P}_{\mathbb H} \Delta: {\mathbb V} \cap \vH^2_{per}(D) \rightarrow {\mathbb H}$.

	
	Let $\pmb{\mathcal{H}}$, $\pmb{\cK}$ be two Hilbert spaces. Then, $\mathcal{L}(\pmb{\mathcal{H}},\pmb{\cK})$ is the space of linear maps from $\pmb{\mathcal{H}}$ to $\pmb{\cK}$. For $m \in \mathbb{N}$, inductively define
	\begin{align}
		\mathcal{L}_m(\pmb{\mathcal{H}}, \pmb{\cK}) := \mathcal{L}(\pmb{\mathcal{H}}, \mathcal{L}_{m-1}(\pmb{\mathcal{H}},\pmb{\cK})),
	\end{align} 
	as the space of all multi-linear maps from $\pmb{\mathcal{H}}\times\cdots\times \pmb{\mathcal{H}}$ ($m$ times) to $\pmb{\cK}$ for $m \geq 2$.
	
	For some function $\vG: \pmb{\mathcal{H}} \rightarrow \pmb{\cK}$, we define the Gateaux derivative of $\vG$ w.r.t. $\vu \in \pmb{\mathcal{H}}$, $D\vG(\vu) \in \mathcal{L}(\pmb{\mathcal{H}},\pmb{\cK})$, whose action is seen as
	\begin{align*}
		\vv \mapsto D\vG(\vu)(\vv)\qquad\forall\vv \in \pmb{\mathcal{H}}.
	\end{align*}
	
	In general, we denote $D^k\vG(\vu)\in \mathcal{L}_m(\pmb{\mathcal{H}},\pmb{\cK})$, as the $k$-Gateaux derivative of $\vG$ w.r.t. $\vu\in \pmb{\mathcal{H}}$.
	
	Throughout this paper, we consider the pairs $\bigl(\pmb{\mathcal{H}},\pmb{\cK}\bigr) = \bigl(\vL^2_{per}(D), \vL^2_{per}(D)\bigr)$, $\bigl(\vH_{per}^1(D),\vH_{per}^1(D)\bigr)$ and $ \bigl(\vH_{per}^2(D),\vH_{per}^2(D)\bigr)$. Below, we state the assumptions on the functional $\vG: \pmb{\mathcal{H}} \rightarrow \pmb{\cK}$.
	\medskip
	\begin{enumerate}
		\item[{\bf(A1)}]  $\vG$ is globally
		Lipschitz continuous and has linear growth. Namely, 
		there exists a constant $C > 0$ such that for all ${\bf v}, {\bf w} \in \pmb{\mathcal{H}}$ 
		\begin{subequations}\label{G}
			\begin{align}\label{Lip}
				\|\vG({\bf v})-\vG({\bf w})\|_{\pmb{\cK}} &\leq C\|{\bf v}-{\bf w}\|_{\pmb{\mathcal{H}}}\, , \\
				\|\vG({\bf v})\|_{\pmb{\cK}}  &\leq C \bigl( \|{\bf v}\|_{\pmb{\mathcal{H}}}+1\bigr)\, .   \label{lineargrow}
			\end{align}
	\end{subequations}
	\item[{\bf (A2)}] There exists a constant $C>0$ such that
	\begin{align}
		\|D\vG\|_{L^{\infty}(\pmb{\mathcal{H}};\mathcal{L}(\pmb{\mathcal{H}},\pmb{\cK}) )} +  \|D^2\vG\|_{L^{\infty}(\pmb{\mathcal{H}};\mathcal{L}_2(\pmb{\mathcal{H}},\pmb{\cK}) )}  \leq C.
	\end{align} 
	\item[{\bf (A3)}] There exists a constant $C>0$ such that for all $\vu,\vv\in \vL_{per}^2(D)$
	\begin{align}
		\|(D\vG(\vu) - D\vG(\vv))\vG(\vv)\|_{\vL^2} \leq C\|\vu - \vv\|_{\vL^2}.
	\end{align}
	\item[{\bf (A4)}] Assume that $\vG: \vH_{per}^1(D) \rightarrow \mV$, i.e. $\vG(\vu)$ is a divergence-free ($\div \vG(\vu) = 0$) for each $\vu \in \vH_{per}^1(D)$.
\end{enumerate}

The assumption ${\bf (A1)}$ is standard for stochastic PDEs. However, the assumptions ${\bf (A2), (A3), (A4)}$ are crucial in proving Lemma \ref{lemma3.1} that is used in obtaining the main theorem of this paper, Theorem \ref{theorem_semi}.

In this paper, we shall use $C$ to denote a generic positive constant
which may depend on $\nu, T$, the datum functions $\vu_0$, $\vf$, and the domain $D$ 
but is independent of the mesh parameter $h$ and $k$. 
In addition, unless it is stated otherwise, for the sake of notation brevity, in the rest of this paper, we set $\vf = 0$ and $\nu = 1$.

\subsection{Variational formulation of problem \eqref{eq1.1}} \label{sec2.1}
We now recall the variational solution concept for \eqref{eq1.1} and refer the reader to \cite{Chow07,PZ1992} for proof of its existence and uniqueness.

\begin{definition}\label{def2.1} 
	Given $(\Omega,\cF, \{\cF_t\},\mP)$, let $W$ be an ${\mathbb R}$-valued Wiener process on it. 
	Suppose ${\bf u}_0\in L^2(\Omega, {\mathbb V})$ and $\vf \in L^2(\Ome;L^2((0,T);L^2_{per}(D)))$.
	An $\{\cF_t\}$-adapted stochastic process  $\{{\bf u}(t) ; 0\leq t\leq T\}$ is called
	a variational solution of \eqref{eq1.1} if ${\bf u} \in  L^2\bigl(\Omega; C([0,T]; {\mathbb V})) 
	\cap L^2\bigl(\Ome;0,T;\vH^2_{per}(D)\bigr)$,
	and satisfies $\mP$-a.s.~for all $t\in (0,T]$
	\begin{align}\label{eq2.8a}
		\bigl({\bf u}(t),  {\bf v} \bigr) + \nu\int_0^t  \bigl(\nab {\bf u}(s), \nab {\bf v} \bigr) 
		\,  ds&=({\bf u}_0, {\bf v})+ \int_0^t \big(\vf(s), \vv\big) \, ds \\\nonumber
		& \qquad+  \bigg(\int_0^t\vG(\vu(s))\,d W(s),\vv\bigg) \qquad\forall  \, {\bf v}\in {\mathbb V}\, . 
	\end{align}

\end{definition}

The following stability estimate for the velocity $\vu$ was proved in \cite{Breit,CP2012}. 

\begin{lemma}\label{stability_pde}
	Let $\vu$ be solution defined in \eqref{eq2.8a}. Assume that $\vu_0 \in L^2\bigl(\Ome; \mV\bigr)$. Then there exists a constant $C = C(\vu_0, \vf, D_T)$ such that
	\begin{enumerate}[\rm (i)]
		\item $\displaystyle 	\mE\biggl[\sup_{0\leq t \leq T} \|\nab\vu(t)\|^2_{\vL^2} + \int_0^T \nu\|\nab^2\vu(t)\|^2_{\vL^2}\, dt \biggr] \leq C.$
		\item  Additionally suppose that $\vu_0\in L^2(\Ome; \mV\cap \vH^2(D))$,
		\begin{align*}
			\mE\biggl[\sup_{0\leq t \leq T} \|\vA\vu(t)\|^2_{\vL^2} + \int_0^T \nu\|\nab^3\vu(t)\|^2_{\vL^2}\, dt \biggr] \leq C.			
		\end{align*}
	\end{enumerate}
	
	
\end{lemma}

\medskip

Next, we introduce the H\"older continuity estimates for the variational solution $\vu$. 

\begin{lemma}\label{lemma2.4}
	Suppose ${\bf u}_0 \in L^2\bigl(\Omega; {\mathbb V\cap \vH^2(D)}\bigr)$. For $\eps >0$, denote $\theta_1 = \frac12 - \epsilon$, $\theta_2 = 1-\epsilon$. There exists a constant $C \equiv C(D_T, \vu_0)>0$, such that the variational solution to problem \eqref{eq1.1} satisfies
	for $s,t \in [0,T]$
	\smallskip
	\begin{enumerate}
		\item[(i)] 	 $\displaystyle{\mathbb E}\bigl[\| {\bf u}(t)-{\bf u}(s)\|^{2}_{\mV} \bigr] 
		\leq C|t-s|^{2\theta_1},$\\
		\item[(ii)] $\displaystyle{\mathbb E}\bigg[\bigg\| {\bf u}(t)-{\bf u}(s) - \int_s^t \vP_{\mH} \vG(\vu(\xi))\, d W(\xi)\bigg\|^{2}_{\mV} \bigg] 
		\leq C|t-s|^{2\theta_2}.$
	\end{enumerate}
\end{lemma}
\begin{proof}
	We just need to give the proof of $(ii)$, and the proof of $(i)$ can be found in \cite{LV2021, CP2012, Breit}. 
	
	Following \cite{LV2021, CP2012, Breit}, we have that the mild solution of \eqref{eq1.1} can be represented as follow.
	\begin{align}
		\vu(t) = e^{-t\vA} \vu_0 + \int_0^t e^{-(t-s)\vA}\vP_{\mH} \vG(\vu(s))\, d W(s).
	\end{align}
	
	For any $t,s \in [0,T]$ with $s<t$, we have
	\begin{align}
		\vu(t) &- \vu(s) - \int_s^t \vP_{\mH}\vG(\vu(\xi))\, dW(\xi) \\\nonumber
		&= \bigl(e^{-t\vA} - e^{-s\vA}\bigr)\vu_0 + \int_0^s \bigl(e^{-(t-\xi)\vA} - e^{-(s-\xi)\vA}\bigr)\vP_{\mH} \vG(\vu(\xi))\, dW(\xi)\\\nonumber
		&\qquad+ \int_{s}^{t} \bigl(e^{-(t-\xi)\vA} - {\bf I}\bigr)\vP_{\mH}\vG(\vu(\xi))\, d W(\xi)\\\nonumber
		&:= {\tt I +II + III}.
	\end{align}
	In order to estimate ${\tt I, II, III}$, we use $\|\cdot\|_{\mV} = \|\vA^{1/2}\cdot\|_{\vL^2}$ and the following inequalities from \cite{PZ1992}:
	\begin{align*}
		\|\vA^a e^{-t\vA}\| \leq C t^{-a}, \qquad \|\vA^{-b}({\bf I} - e^{-t\vA})\| \leq C t^b.
	\end{align*}
	
	First, we have
	\begin{align}
		\|{\tt I}\|_{\mV} &= \|e^{-s\vA}\vA^{\theta_2} \vA^{-\theta_2}\bigl(e^{-(t-s)\vA} - {\bf I}\bigr)\vA^{1/2}\vu_0\|_{\vL^2}\\\nonumber
		&\leq Cs^{-\theta_2} (t-s)^{\theta_2}\|\vA^{1/2} \vu_0\|_{\vL^2}.
	\end{align}
	So, 
	\begin{align}
		\mE\bigl[\|{\tt I}\|^2_{\mV}\bigr] \leq C(t-s)^{2\theta_2}\mE\bigl[\|\vu_0\|^2_{\mV}\bigr].
	\end{align}	
	
	By using the It\^o isometry and assumptions on $\vG$, we obtain
	\begin{align}
		\mE\bigl[\|{\tt II}\|^2_{\mV}\bigr] &= \mE\bigg[\int_0^s\|\vA^{\frac12}\bigl(e^{-(t-\xi)\vA} - e^{-(s-\xi)\vA}\bigr)\vP_{\mH}\vG(\vu(\xi))\|^2_{\vL^2}\, d\xi\bigg]\\\nonumber
		&=\mE\bigg[\int_0^s\|\vA^{\frac12}e^{-(s-\xi)\vA}\vA^{\theta_2}\vA^{-\theta_2}\bigl(e^{-(t- s)\vA} - {\bf I}\bigr)\vP_{\mH}\vG(\vu(\xi))\|^2_{\vL^2}\, d\xi\bigg]\\\nonumber
		&=\mE\bigg[\int_0^s\|e^{-(s-\xi)\vA}\vA^{\frac12-\eps}\vA^{-\theta_2}\bigl(e^{-(t- s)\vA} - {\bf I}\bigr)\vA\vP_{\mH}\vG(\vu(\xi))\|^2_{\vL^2}\, d\xi\bigg]\\\nonumber
		&\leq C(t-s)^{2\theta_2} \mE\bigl[\sup_{\xi\in [0,T]}\|\vu(\xi)\|^2_{\vH^2}\bigr] \int_0^s \frac{d\xi}{(s-\xi)^{1-2\eps}}\\\nonumber
		&\leq C(t-s)^{2\theta_2}\mE\bigl[\sup_{\xi\in [0,T]}\|\vu(\xi)\|^2_{\vH^2}\bigr].
	\end{align}
	
	In order to estimate {\tt III}, we also use the It\^o isometry and assumtions on $\vG$, we have
	\begin{align}
		\mE\bigl[\|{\tt III}\|^2_{\mV}\bigr] &= \int_{s}^t \mE\bigr[\|\vA^{\frac12}\bigl(e^{-(t-\xi)\vA} - {\bf I}\bigr)\vP_{\mH}\vG(\vu(\xi))\|^2_{\vL^2}\bigr]\, d\xi\\\nonumber
		&= \int_{s}^t \mE\bigr[\|\vA^{-\frac12}\bigl(e^{-(t-\xi)\vA} - {\bf I}\bigr)\vA\vP_{\mH}G(\vu(\xi))\|^2_{\vL^2}\bigr]\, d\xi\\\nonumber
		&\leq C\int_{s}^t (t-\xi) \mE\bigl[\|\vu(\xi)\|^2_{\vH^2}\bigr]\, d\xi\\\nonumber
		&\leq C(t-s)^2\mE\bigl[\sup_{\xi\in [0,T]}\|\vu(\xi)\|^2_{\vH^2}\bigr].
	\end{align}
	
	The desired estimate is obtained by combining all the estimates of {\tt I, II, III}. The proof is complete.
	
\end{proof}

\begin{remark} Lemma \ref{lemma2.4} $(ii)$ provides a new H\"older continuity estimate for the velocity solution of \eqref{eq1.1}. This estimate plays a very important role in the proof of Theorem \ref{theorem_semi} which is the main result of this paper.  
\end{remark}
\subsection{Interpretation of the pressure solution}\label{sec2.2}
Definition \ref{def2.1} only defines the velocity $\mathbf{u}$ for \eqref{eq1.1}, 
its associated pressure $p$, in contrast, is very low regularity. In fact, according to \cite[Theorem 2.2]{LRS03}, the pressure solution $p \in L^1(\Omega; W^{-1.\infty}(0,T; L^2(D))$ is a distribution that is still very challenging for obtaining an error estimate of any numerical scheme.  In that regard, we quote the following 
theorem from \cite{FPL2021,FL2020} that comes with a new interpretation for the pressure solution.

\begin{theorem}\label{thm 2.2}
	Let $\{{\bf u}(t) ; 0\leq t\leq T\}$ be the variational solution of \eqref{eq1.1}. There exists a unique adapted process 
	{$P\in {L^2\bigl(\Omega; L^2(0,T; H^1_{per}(D)/{\mathbb R})\bigr)}$} such that $(\mathbf{u}, P)$ satisfies 
	$\mP$-a.s.~for all $t\in (0,T]$
	\begin{subequations}\label{eq2.100}
		\begin{align}\label{eq2.10a}
			&\bigl({\bf u}(t),  \pphi \bigr) + \nu\int_0^t  \bigl(\nab {\bf u}(s), \nab \pphi \bigr) \, ds
			- \bigl(  \div \pphi, P(t) \bigr) \\
			&=({\bf u}_0, \pphi) + \int_0^t \big(\vf(s), \pphi\big) \, ds 
			+  \bigg(\int_0^t\vG(\vu(s))\, dW(s),\pphi\bigg)  \,\,\, \forall  \, \pphi\in \vH^1_{per}(D)\, , \nonumber \\ 
			&\bigl(\div {\bf u}, q \bigr) =0 \qquad\forall \, q\in L^2_{per}(D)/\mathbb{R}   .  \label{eq2.10b}
		\end{align}
	\end{subequations}
\end{theorem}

Then, the pressure solution $p$ is defined as the time distributional derivative of the ``time-averaged" pressure $P$. To the best of our knowledge, all the pressure error estimates in the literature \cite{FPL2021,FL2020,FQ} until now are obtained for $P$ only. However, these error estimates just provide a weak convergence. Moreover, it is also not very easy to simulate the pressure solution $P$ in numerical tests since $P$ must be approximated by a high-order quadrature approximation. Therefore, the second goal of this paper is to introduce a ``simulate" pressure solution $\tilde{p}$ that is defined as follows:

For $\delta>0 $ such that $t-\delta\geq 0$ for all $t \in (0,T]$, define
\begin{align}\label{simu_press}
	\tilde{p}(t) := \frac{P(t) - P(t - \delta)}{\delta} \qquad\forall t \in (0,T].
\end{align}

With this ``simulate" pressure solution, we can rewite \eqref{eq2.100} as follows
\begin{subequations}
	\begin{align}
		&\bigl({\bf u}(t) - \vu(t-\delta),  \pphi \bigr) + \nu\int_{t-\delta}^{t}  \bigl(\nab {\bf u}(s), \nab \pphi \bigr) \, ds
		- \delta \bigl(  \div \pphi, \tilde{p}(t) \bigr) \\
		&= \int_{t-\delta}^{t} \big(\vf(s), \pphi\big) \, ds 
		+  \bigg(\int_{t-\delta}^{t}\vG(\vu(s))\, d W(s),\pphi\bigg)  \,\,\, \forall  \, \pphi\in \vH^1_{per}(D)\, , \nonumber \\ 
		&\bigl(\div {\bf u}, q \bigr) =0 \qquad\forall \, q\in L^2_{per}(D)/\mathbb{R}   .  
	\end{align}
\end{subequations}

\section{Semi-discretization in time}\label{section_semi} In this section, we follow the strategy of the Milstein scheme in SDEs to propose a new time discretization method of \eqref{eq2.8a}. 
\subsection{Formulation of the proposed method}

Let $I_k: = \{t_n\}_{n=1}^M$ be a uniform mesh of the interval $[0,T]$ 
with the time step-size $k = \frac{T}{M}$. Note that $t_0 = 0$ and $t_M = T$. 

\smallskip
\noindent
\textbf{Algorithm 1} 

Let $\vu^0 = \vu_0$ be a given $\mV$-valued random variable. Find the pair $\{\vu^{n+1}, p^{n+1}\} \in \mV\times L^2_{per}(D)/\mathbb{R}$ recursively such that $\mP$-a.s. 
\begin{subequations}\label{milsteinscheme}
	\begin{align}
		\label{milsteinscheme1}	\bigl(\vu^{n+1} - \vu^n, \pphi\bigr) &+ \nu k \bigl(\nab\vu^{n+1}, \nab\pphi\bigr) - k\bigl(p^{n+1}, \div \pphi\bigr) \\\nonumber &= \bigl(\vG(\vu^n)\Delta W_{n} + \frac12 D\vG(\vu^n)\,\vG(\vu^n)\bigl[(\Delta W_n)^2 - k\bigr],\pphi\bigr),\\
		\label{milsteinscheme2}	\bigl(\div \vu^{n+1},\psi\bigr) &=0, 
	\end{align}
\end{subequations}
for all $\pphi \in \vH^1_{per}(D)$ and $\psi \in L^2_{per}(D)$ and $\Delta W_n = W(t_{n+1}) - W(t_n) \sim \mathcal{N}(0,k)$.


Next, we define $\displaystyle \G: \mathbb{R}^+\times \vH_{per}^1(D) \rightarrow \vL^2_{per}(D)$ by
\begin{align}\label{eq3.2}
	\G(s;\vu) := \vG(\vu) + D\vG(\vu)\vG(\vu)\, \int_{t_{n}}^s\,d W(r),\qquad t_n \leq s \leq t_{n+1}
\end{align}

Then we have
\begin{align*}
	\int_{t_{n}}^{t_{n+1}}\G(s;\vu^n)\, dW(s) &= \vG(\vu^n)\Delta W_{n} + D\vG(\vu^n)\vG(\vu^n) \int_{t_{n}}^{t_{n+1}}\int_{t_{n}}^s d W(r)\, dW(s)\\\nonumber
	&= \vG(\vu^n)\Delta W_{n} + \frac12D\vG(\vu^n)\vG(\vu^n) \bigl[(\Delta W_{n})^2 - k\bigr].
\end{align*}
Therefore, we rewrite \eqref{milsteinscheme1} as follow:
\begin{align}\label{milsteinscheme3}
	\bigl(\vu^{n+1} - \vu^n, \pphi\bigr) + \nu k \bigl(\nab\vu^{n+1}, \nab\pphi\bigr) &- k\bigl(p^{n+1}, \div \pphi\bigr)  \\\nonumber&= \bigg(\int_{t_{n}}^{t_{n+1}}\G(s;\vu^n)\, dW(s) ,\pphi\bigg),
\end{align}

\medskip

Next, we state the following technical lemma that is used to prove the error estimate of the velocity approximation in Theorem \ref{theorem_semi}.

\begin{lemma}\label{lemma3.1} Suppose that $\vG$ satisfies the assumptions ${\bf(A1), (A2), (A3), (A4)}$. Let $\vu_0\in L^2(\Ome;\mV\cap \vH^2(D))$, there exist constants $C>0$ such that the functional $\G$ defined in \eqref{eq3.2} satisfies 
	\smallskip
	\begin{enumerate}[\rm (i)]
		\item  $\displaystyle	\|\G(s;\vu) - \G(s;\vv)\|_{\vL^2} \leq C\|\vu - \vv\|_{\vL^2},\qquad\forall s>0,\vu,\vv \in \vL_{per}^2(D)$,
		\smallskip
		\item  $\displaystyle \mE\bigl[\bigl\|\vG(\vu(s)) - \G(s;\vu(t_n))\bigr\|^2_{\vL^2}\bigr] \leq C|s - t_n|^{2(1-\epsilon)}$, for $t_n\leq s < t_{n+1}$ and $\epsilon>0$.
	\end{enumerate}
\end{lemma}

\smallskip
\begin{proof}
	The Lipschitz continuity of $\mathcal{\G}$ in $(i)$ is directly obtained from \eqref{Lip} and assumptions of $\vG$. Below, we give the proof of $(ii)$.
	First, let's denote 
	\begin{align}\label{eq3.4}
		{\bf R_1} : &= \vu(s) - \vu(t_n) - \int_{t_{n}}^s \vG(\vu(\xi))\, dW(\xi).
	\end{align}
	Then by using Lemma \ref{lemma2.4} $(ii)$, we obtain
	\begin{align}\label{eq3.55}
		\mE\bigl[\|{\bf R_1}\|^2_{\vL^2}\bigr] \leq C|t_n - s|^{2(1-\eps)}.
	\end{align}
	By using the Taylor expansion for $\vG$, we have
	\begin{align}\label{eq3.5}
		\vG(\vu(s)) = \vG(\vu(t_n)) + D\vG(\vu(t_n))(\vu(s) - \vu(t_n)) + {\bf R_2},
	\end{align}
	where  \begin{align*}
		{\bf R_2} := \int_0^1 (1-\eta) \bigl(D^2\vG(\vu(t_n) + \eta(\vu(s) - \vu(t_n)))(\vu(s) - \vu(t_n))\bigr)(\vu(s) - \vu(t_n))\, d\eta.
	\end{align*}
	Substituting \eqref{eq3.4} into \eqref{eq3.5}, we have
	\begin{align}\label{eq3.6}
		\vG(\vu(s)) &= \vG(\vu(t_n)) + D\vG(\vu(t_n))\Bigl({\bf R_1} +  \int_{t_{n}}^s \vG(\vu(\xi))\, dW(\xi)\Bigr) + {\bf R_2}\\\nonumber
		&=\vG(\vu(t_n)) + D\vG(\vu(t_n)) \int_{t_{n}}^s \vG(\vu(\xi))\, dW(\xi)+ D\vG(\vu(t_n)){\bf R_1} + {\bf R_2}.
	\end{align}
	Next, by \eqref{eq3.2}, we have
	\begin{align}\label{eq3.77}
		\G(s;\vu(t_n)) = \vG(\vu(t_n)) + D\vG(\vu(t_n))\vG(\vu(t_n))\bigl(W(s) - W(t_n)\bigr).
	\end{align}
	From \eqref{eq3.6} and \eqref{eq3.77} we obtain
	\begin{align}
		\vG(\vu(s)) - \G(s;\vu(t_n)) &= D\vG(\vu(t_n)) \int_{t_{n}}^s \bigl[\vG(\vu(\xi)) - \vG(\vu(t_n))\bigr]\, dW(\xi)\\\nonumber
		&\qquad+ D\vG(\vu(t_n)){\bf R_1} +  {\bf R_2}.
	\end{align}
	Therefore, by using the assumption ${\bf (A2)}$, and the It\^o isometry we have
	\begin{align}\label{eq3.99}
		\mE\bigr[\|\vG(\vu(s)) - \G(s;\vu(t_n))\|^2_{\vL^2}\bigr] &\leq C\int_{t_{n}}^s\mE\bigl[\|\vG(\vu(\xi)) - \vG(\vu(t_n))\|^2_{\vL^2}\bigr]\, d\xi \\\nonumber
		&\qquad+ C\mE\bigl[\|{\bf R_1}\|^2_{\vL^2}\bigr] + \mE\bigl[\|{\bf R_2}\|^2_{\vL^2}\bigr].
	\end{align}
	The last term on the right side of \eqref{eq3.99} is controlled by \eqref{eq3.55}. In addition, by ${\bf (A1)}$ and Lemma \ref{lemma2.4} $(i)$, we have
	\begin{align}
		C\int_{t_{n}}^s\mE\bigl[\|\vG(\vu(\xi)) - \vG(\vu(t_n))\|^2_{\vL^2}\bigr]\, d\xi &\leq C\int_{t_{n}}^s \mE\bigl[\|\vu(\xi) - \vu(t_n)\|^2_{\vL^2}\bigr] \, d\xi\\\nonumber
		&\leq C|t_n - s|^{2(1-\eps) }.
	\end{align}
	Finally, by using again ${\bf (A2)}$ and Lemma \ref{lemma2.4} $(i)$ we obtain
	\begin{align*}
		\mE[\|{\bf R_2}\|^2_{\vL^2}] \leq C\mE\bigl[\|\vu(s) - \vu(t_n)\|^4_{\vL^4}\bigr] \leq C|t_n-s|^{2(1-\eps)}.
	\end{align*}
	The proof is complete.

\end{proof}

\medskip

Next, we state the following stability estimates for the velocity approximation $\{\vu^n\}$ of Algorithm 1.

\begin{lemma}\label{stability_mean}
	Let $\vu_0 \in L^{2}(\Ome;\mV)$ and $\vG$ satisfies ${\bf (A1), (A2),  (A4)}$. Then there exists a constant $C = C(\vu_0,D_T)$ such that the following estimations hold:
	\begin{enumerate}[\rm (i)]
		\item\qquad $\displaystyle\mE\biggl[\max_{1\leq n \leq M}\|\vu^n\|^{2}_{\mV} + \nu k\sum_{n=1}^M \|{\bf A}\vu^n\|^2_{\vL^2}\biggr] \leq C$.
		\item 	\qquad$\displaystyle\mE\biggl[k\sum_{n=1}^M \|\nab p^{n}\|^2_{\vL^2}\biggr] \leq C$.
	\end{enumerate}
\end{lemma}

\begin{proof}
	The proof of Lemma \ref{stability_mean} is similar to the proof of \cite[Lemma 3.1]{FPL2021, CP2012}. We omit the detailed proof to save space.
\end{proof}

\subsection{Error estimates for Algorithm 1}\label{sec3.2}

In this part, we state the main result of this paper which establishes an $O(k^{1-\eps})$ convergence order for the velocity approximation of the Milstein method.

\begin{theorem}\label{theorem_semi} 
	Let $\vu$ be the variational solution to \eqref{eq2.8a} and $\{\vu^{n}\}_{n=1}^M$ be generated by  Algorithm 1. Assume that $\vG$ satisifies ${\bf (A1), (A2), (A3), (A4)}$ and $\vu_0 \in L^{2}(\Ome; \mV\cap \vH^2(D))$. For $\epsilon>0$, there exists $C = C(\vu_0,D_T)>0$ such that 
	\begin{align}\label{eq310}
		\Bigl(\mE\Bigl[\max_{1\leq n \leq M} \|\vu(t_n) - \vu^n\|^2_{\vL^2}\Bigr]\Bigr)^{\frac12} + \bigg(\mE\bigg[k\sum_{n=1}^M\|\nab(\vu(t_n) - \vu^n)\|^2_{\vL^2}\bigg]\bigg)^{\frac12} \leq C\, k^{1-\epsilon}
	\end{align}
\end{theorem}

\begin{proof}
	Denote $\ve^{n} := \vu(t_n) - \vu^n$, then by subtracting \eqref{eq2.8a} from \eqref{milsteinscheme3}, we obtain
	\begin{align}\label{err_eq}
		\bigl(\ve^{n+1}- \ve^n, \pphi\bigr) + k\bigl(\nab\ve^{n+1},\nab\pphi\bigr) &= \int_{t_{n}}^{t_{n+1}}\bigl(\nab(\vu(t_{n+1}) - \vu(s)),\nab\pphi\bigr)\, ds\\\nonumber
		&+ \bigg(\int_{t_{n}}^{t_{n+1}}(\vG(\vu(s)) - \G(s;\vu^n))\, dW(s),\pphi\bigg) . 
	\end{align}
	Choose $\pphi = \ve^{n+1} \in \mV$ in \eqref{err_eq} and use the identity $2a(a-b) = a^2 - b^2 + (a-b)^2$, we obtain
	\begin{align}\label{eq3.8}
		\frac12\bigl[\|\ve^{n+1}\|^2_{\vL^2} - \|\ve^n\|^2_{\vL^2}\bigr] &+ \frac12\|\ve^{n+1} - \ve^n\|^2_{\vL^2} + k\|\nab\ve^{n+1}\|^2_{\vL^2}\\\nonumber
		&=\int_{t_{n}}^{t_{n+1}}\bigl(\nab(\vu(t_{n+1}) - \vu(s)),\nab\ve^{n+1}\bigr)\, ds\\\nonumber
		&\qquad+ \bigg(\int_{t_{n}}^{t_{n+1}}(\vG(\vu(s)) - \G(s;\vu^n))\, dW(s),\ve^{n+1}\bigg)  \\\nonumber
		&:= {\tt I + II}.
	\end{align}
	Next, we estimate the right-hand side of \eqref{eq3.8} as follows.
	\begin{align}\label{eq3.9}
		{\tt I} &= \int_{t_{n}}^{t_{n+1}}\bigg(\nab\bigg(\vu(t_{n+1}) - \vu(s) - \int_s^{t_{n+1}} \vG(\vu(\xi))\, dW(\xi)\bigg),\nab\ve^{n+1}\bigg)\, ds \\\nonumber
		&\qquad+ \int_{t_n}^{t_{n+1}} \bigg(\int_s^{t_{n+1}}\nab \vG(\vu(\xi))\, dW(\xi),\nab\ve^{n+1}\bigg)\, ds\\\nonumber
		&= \int_{t_{n}}^{t_{n+1}}\bigg(\nab\bigg(\vu(t_{n+1}) - \vu(s) - \int_s^{t_{n+1}} \vG(\vu(\xi))\, dW(\xi)\bigg),\nab\ve^{n+1}\bigg)\, ds \\\nonumber
		&\qquad+ \int_{t_n}^{t_{n+1}} \bigg(\int_s^{t_{n+1}}\nab \vG(\vu(\xi))\, dW(\xi),\nab(\ve^{n+1} - \ve^n)\bigg)\, ds \\\nonumber
		&\qquad+ \int_{t_n}^{t_{n+1}} \bigg(\int_s^{t_{n+1}}\nab \vG(\vu(\xi))\, dW(\xi),\nab\ve^n\bigg)\, ds\\\nonumber
		&:= {\tt I_1 + I_2 + I_3}.
	\end{align}
	Now, we bound the right-hand side of \eqref{eq3.9} as follows:
	
	By using Cauchy-Schwarz inequality and then Lemma \ref{lemma2.4}$(ii)$,we obtain
	\begin{align}\label{eq3.10}
		\mE[{\tt I_1}] &\leq \mE\bigg[\int_{t_n}^{t_{n+1}} \bigg\|\vu(t_{n+1}) -\vu(s) - \int_s^{t_{n+1}}\vG(\vu(\xi))\, dW(\xi)\bigg\|_{\mV}^2\, ds\bigg] \\\nonumber
		&\qquad\qquad\qquad\qquad\qquad\qquad+ \frac{k}{4}\mE\bigl[\|\nab\ve^{n+1}\|^2_{\vL^2}\bigr]\\\nonumber
		&\leq Ck^{1+2(1-\epsilon)} + \frac{k}{4}\mE\bigl[\|\nab\ve^{n+1}\|^2_{\vL^2}\bigr].
	\end{align}
	
	By doing integration by parts and applying H\"older inequality, we obtain
	\begin{align}
		{\tt I_2} &= -\int_{t_n}^{t_{n+1}}\bigg(\int_s^{t_{n+1}}\Delta\vG(\vu(\xi))\, dW(\xi), \ve^{n+1} - \ve^n\bigg)\, ds\\\nonumber
		&\leq 2\bigg\|\int_{t_n}^{t_{n+1}}\int_s^{t_{n+1}}\Delta G(\vu(\xi))\,dW(\xi)\,ds\bigg\|^2_{\vL^2} + \frac18\|\ve^{n+1} - \ve^n\|^2_{\vL^2}\\\nonumber
		& = 2\int_D\bigg|\int_{t_n}^{t_{n+1}}\int_s^{t_{n+1}}\Delta \vG(\vu(\xi))\,dW(\xi)\,ds\bigg|^2\, d{\bf x}+ \frac18\|\ve^{n+1} - \ve^n\|^2_{\vL^2}\\\nonumber
		&\leq 2\int_D\bigg(\int_{t_n}^{t_{n+1}}\bigg|\int_s^{t_{n+1}}\Delta\vG(\vu(\xi))\,dW(\xi)\bigg|\,ds\bigg)^2\, d{\bf x}+ \frac18\|\ve^{n+1} - \ve^n\|^2_{\vL^2}\\\nonumber
		&\leq 2k\int_D\int_{t_n}^{t_{n+1}}\bigg|\int_s^{t_{n+1}}\Delta\vG(\vu(\xi))\, dW(\xi)\bigg|^2\, ds\, d{\bf x}  + \frac18\|\ve^{n+1} - \ve^n\|^2_{\vL^2}\\\nonumber
		&= 2k \int_{t_n}^{t_{n+1}}\bigg\|\int_s^{t_{n+1}}\Delta\vG(\vu(\xi))\, dW(\xi)\bigg\|^2_{\vL^2}\, ds + \frac18\|\ve^{n+1} - \ve^n\|^2_{\vL^2}.
	\end{align}
	
	So, by using the It\^o isometry, we have
	\begin{align}
		\mE[{\tt I_2}] &\leq 2k\mE\bigg[\int_{t_n}^{t_{n+1}}\int_s^{t_{n+1}}\|\Delta\vG(\vu(\xi))\|^2_{\vL^2}\,d\xi\,ds\bigg] + \frac18\mE\bigl[\|\ve^{n+1} - \ve^n\|^2_{\vL^2}\bigr]\\\nonumber
		&\leq Ck^3\mE\bigl[\sup_{\xi\in [0,T]} \|\vu(\xi)\|^2_{\vH^2}\bigr] + \frac18\mE\bigl[\|\ve^{n+1} - \ve^n\|^2_{\vL^2}\bigr].
	\end{align}
	
	By using the martingale property of the It\^o integral, we conclude that $\mE[{\tt I_3}] = 0$. Finally, from the estimates of ${\tt I_1, I_2, I_3}$, we obtain an estimate for the term {\tt I} as below
	\begin{align}
		\mE[{\tt I}] \leq Ck^{1+2(1-\epsilon)} + \frac{k}{4}\mE\bigl[\|\nab\ve^{n+1}\|^2_{\vL^2}\bigr] &+ Ck^3\mE\bigl[\sup_{\xi\in [0,T]} \|\vu(\xi)\|^2_{\vH^2}\bigr] \\\nonumber
		&+ \frac18\mE\bigl[\|\ve^{n+1} - \ve^n\|^2_{\vL^2}\bigr].
	\end{align}
	
	Next, we estimate {\tt II} as follows. 
	\begin{align}
		{\tt II} &= \bigg(\int_{t_{n}}^{t_{n+1}}(\vG(\vu(s)) - \G(s;\vu^n))\, dW(s),\ve^{n+1} - \ve^{n}\bigg) \\\nonumber
		&\qquad+ \bigg(\int_{t_{n}}^{t_{n+1}}(\vG(\vu(s)) - \G(s;\vu^n))\, dW(s),\ve^{n}\bigg)  \\\nonumber
		&:= {\tt II_1 + II_2}.
	\end{align}
	
	First, we observe that $\mE[{\tt II_2}] =0$ due to the martingale property of the It\^o integral. So, we just need to estimate ${\tt II_1}$. Bu using again the Cauchy-Schwarz inequality we obtain
	\begin{align}
		{\tt II_1} &\leq 2\bigg\|\int_{t_n}^{t_{n+1}} \bigl(\vG(\vu(s)) - \G(s;\vu^n)\bigr)\, dW(s)\bigg\|^2_{\vL^2} + \frac18\|\ve^{n+1} - \ve^n\|^2_{\vL^2} .
	\end{align}
	
	Taking the expectation and applying the It\^o isometry, we obtain
	\begin{align}\label{eq3.16}
		\mE[{\tt II_1}] &\leq 2\mE\bigg[\int_{t_n}^{t_{n+1}} \|\vG(\vu(s)) - \G(s;\vu^n)\|^2_{\vL^2}\, ds\bigg] + \frac18\mE\bigl[\|\ve^{n+1} - \ve^n\|^2_{\vL^2}\bigr]\\\nonumber
		&\leq 2\mE\bigg[\int_{t_n}^{t_{n+1}} \|\vG(\vu(s)) - \G(s;\vu(t_n))\|^2_{\vL^2}\, ds\bigg] \\\nonumber
		&\quad+  2\mE\bigg[\int_{t_n}^{t_{n+1}} \|\vG(\vu(t_n)) - \G(s;\vu^n)\|^2_{\vL^2}\, ds\bigg] + \frac18\mE\bigl[\|\ve^{n+1} - \ve^n\|^2_{\vL^2}\bigr]\\\nonumber
		&\leq Ck^{1 + 2(1-\epsilon)} + Ck\mE\bigl[\|\ve^n\|^2_{\vL^2}\bigr] + \frac18\mE\bigl[\|\ve^{n+1} - \ve^n\|^2_{\vL^2}\bigr],
	\end{align}
	where the last inequality of \eqref{eq3.16} was obtained by using Lemma \ref{lemma3.1}  $(i), (ii)$.
	
	Substituting the estimates of ${\tt I, II}$ into \eqref{eq3.8} in the expectation we obtain
	\begin{align}\label{eq3.17}
		\frac12\mE\bigl[&\|\ve^{n+1}\|^2_{\vL^2} - \|\ve^n\|^2_{\vL^2}\bigr] + \frac14\mE\bigl[\|\ve^{n+1} - \ve^n\|^2_{\vL^2}\bigr] + \frac{3k}{4}\mE\bigl[\|\nab\ve^{n+1}\|^2_{\vL^2}\bigr] \\\nonumber
		&\leq Ck^{1+2(1-\eps)} + Ck^3\mE\bigl[\sup_{\xi\in [0,T]}\|\vu(\xi)\|^2_{\vH^2}\bigr] + Ck\mE\bigl[\|\ve^n\|^2_{\vL^2}\bigr].
	\end{align} 
	
	Next, taking the summation $\sum_{n=0}^m$, for any $0\leq m < M$, we obtain
	\begin{align}\label{eq3.18}
		\frac12\mE\bigl[&\|\ve^{m+1}\|^2_{\vL^2}\bigr] + \frac14\sum_{n=0}^m\mE\bigl[\|\ve^{n+1} - \ve^n\|^2_{\vL^2}\bigr] + \frac{3}{4}k\sum_{n=0}^m\mE\bigl[\|\nab\ve^{n+1}\|^2_{\vL^2}\bigr] \\\nonumber
		&\leq Ck^{2(1-\eps)} + Ck^2\mE\bigl[\sup_{\xi\in [0,T]}\|\vu(\xi)\|^2_{\vH^2}\bigr] + Ck\sum_{n=0}^m\mE\bigl[\|\ve^n\|^2_{\vL^2}\bigr].
	\end{align} 
	
	By applying the discrete Grownwall inequality to \eqref{eq3.18}, we obtain
	\begin{align}\label{eq3.19}
		\frac12\mE\bigl[&\|\ve^{m+1}\|^2_{\vL^2}\bigr] + \frac14\sum_{n=0}^m\mE\bigl[\|\ve^{n+1} - \ve^n\|^2_{\vL^2}\bigr] + \frac{3}{4}k\sum_{n=0}^m\mE\bigl[\|\nab\ve^{n+1}\|^2_{\vL^2}\bigr] \\\nonumber 
		&\leq Ck^{2(1-\eps)}\bigl(1 + k^{2\eps}\mE\bigl[\sup_{\xi\in [0,T]} \|\vu(\xi)\|^2_{\vH^2}\bigr]\bigr)\exp(Ct_m).
	\end{align}
	Next, taking $\max_{0\leq n \leq M-1}$ on \eqref{eq3.19} we obtain
	\begin{align}\label{eq3.26}
		\max_{1\leq m \leq M} \mE\bigl[\|\ve^m\|^2_{\vL^2}\bigr] + \mE\Bigl[k\sum_{n = 1}^M \|\nab\ve^m\|^2_{\vL^2}\Bigr] \leq Ck^{2(1-\eps)}.
	\end{align}
	Next, we will use \eqref{eq3.26} to establish the desired estimate \eqref{eq310}. To do that, we start with applying the summation $\sum_{n = 1}^m$, for some $1 \leq m \leq M$, to \eqref{eq3.8} and then take the maximum overall $m's$ as well as the expectation. These processes yield to
	\begin{align}\label{eq3.27}
		&\frac12\mE\Bigl[\max_{1\leq m \leq M}\|\ve^m\|^2_{\vL^2}\Bigr] + \frac12\mE\Bigl[\sum_{n = 1}^M \|\ve^n - \ve^{n-1}\|^2_{\vL^2}\Bigr] + \mE\Bigl[k\sum_{n = 1}^M \|\nab\ve^n\|^2_{\vL^2}\Bigr]\\\nonumber
		&\leq\mE\Bigl[\max_{1\leq m \leq M}\Bigl|\sum_{n = 1}^m \int_{t_{n-1}}^{t_{n}}\bigl(\nab(\vu(t_{n}) - \vu(s)),\nab\ve^{n}\bigr)\, ds\Bigr|\Bigr]\\\nonumber
		&\qquad+\mE\bigg[\max_{1\leq m \leq M}\bigg|\sum_{n = 1}^m \bigg(\int_{t_{n-1}}^{t_{n}}(\vG(\vu(s)) - \G(s;\vu^{n-1}))\, dW(s),\ve^{n}\bigg) \bigg|\bigg]\\\nonumber
		&:=A + B
	\end{align} 
	
	Now, we estimate the right-hand side of \eqref{eq3.27} as follows. To estimate $A$, we follow the same strategy shown in estimating ${\tt I}$ in \eqref{eq3.9}. Namely, by using \eqref{eq3.26} we obtain
	\begin{align*}
		A &\leq \mE\bigg[\sum_{n = 1}^M\int_{t_{n-1}}^{t_{n}}\bigg\|\nab\Bigl(\vu(t_{n}) - \vu(s) - \int_s^{t_{n}} \vG(\vu(\xi))\, dW(\xi)\Bigr)\bigg\|_{\vL^2}\|\nab\ve^{n}\|_{\vL^2}\, ds\bigg] \\\nonumber
		&\qquad+ \mE\Bigl[\max_{1\leq m \leq M}\bigg|\sum_{n = 1}^m\int_{t_{n-1}}^{t_{n}} \bigg(\int_s^{t_{n}}\Delta \vG(\vu(\xi))\, dW(\xi),\ve^{n} - \ve^{n-1}\bigg)\, ds\bigg|\bigg] \\\nonumber
		&\qquad+\mE\bigg[ \max_{1\leq n \leq M}\bigg|\sum_{n=1}^m\int_{t_{n-1}}^{t_{n}} \bigg(\int_s^{t_{n}}\Delta\vG(\vu(\xi))\, dW(\xi),\ve^{n-1}\bigg)\, ds\bigg|\bigg]. 
	\end{align*}
	
	Next, by using Fubini's theorem and the Burholder-Davis-Gundy inequality \cite[Lemma 4.3]{Chow07} as well as the estimate \eqref{eq3.26} and Lemma \ref{lemma2.4} $(ii)$ we obtain
	\begin{align}
		A &\leq \mE\bigg[\sum_{n = 1}^M \int_{t_{n-1}}^{t_n} \Bigl\|\vu(t_n) - \vu(s) - \int_{s}^{t_n} \vG(\vu(\xi))\, dW(\xi)\Bigr\|^2_{\mV} + \frac14 k\sum_{n = 1}^M \|\nab\ve^n\|^2_{\vL^2}\bigg] \\\nonumber
		&+ \mE\bigg[\sum_{n = 1}^M \Bigl(\Bigl\|\int_{t_{n-1}}^{t_{n}}\int_{s}^{t_n} \Delta \vG(\vu(\xi))\, dW(\xi)\, ds\Bigr\|^2_{\vL^2} + \frac{1}{4}\|\ve^n - \ve^{n-1}\|^2_{\vL^2}\Bigr)\bigg]\\\nonumber
		&+ \mE\bigg[ \max_{1\leq n \leq M}\bigg|\sum_{n=1}^m \bigg(\int_{t_{n-1}}^{t_{n}}\int_{t_{n-1}}^{\xi}\Delta\vG(\vu(\xi))\,ds\,dW(\xi) ,\ve^{n-1}\bigg)\bigg|\bigg]\\\nonumber
		&:= {\tt a + b +c}.
	\end{align}
	By using Lemma \ref{lemma2.4} $(ii)$ and the estimate \eqref{eq3.26} we obtain
	\begin{align*}
		{\tt a + b} \leq Ck^{2(1-\eps)} + Ck^{2(1-\eps)}\sup_{\xi\in [0,T]}\mE\bigl[\|\vu(\xi)\|^2_{\vH^2}\bigr] \leq Ck^{2(1-\eps)}.
	\end{align*}
	Next, by using the Burholder-Davis-Gundy inequality \cite[Lemm 4.3]{Chow07} and \eqref{eq3.26} we have
	\begin{align}
		{\tt c} 	&\leq \mE\bigg[\Bigl(\sum_{n = 1}^M\int_{t_{n-1}}^{t_n} \Bigl\|\int_{t_{n-1}}^{\xi}\Delta \vG(\vu(\xi))\, ds\Bigr\|^2_{\vL^2}\|\ve^{n-1}\|^2_{\vL^2}\, d\xi\Bigr)^{1/2}\bigg]\\\nonumber
		&\leq \mE\bigg[\frac18\max_{1\leq n \leq M}\|\ve^n\|^2_{\vL^2} + 2\sum_{n = 1}^M \int_{t_{n-1}}^{t_n} \Bigl\|\int_{t_{n-1}}^{\xi}\Delta \vG(\vu(\xi))\, ds\Bigr\|^2_{\vL^2}\, d\xi\bigg]\\\nonumber
		&\leq \frac18\mE\bigl[\max_{1\leq n \leq M}\|\ve^n\|^2_{\vL^2}\bigr] + Ck^{2(1-\eps)}\mE\bigl[\sup_{\xi\in [0,T]}\|\vu(\xi)\|^2_{\vH^2}\bigr].
	\end{align}
	By using the Burholder-Davis-Gundy inequality and Lemma \ref{lemma3.1}, we obtain
	\begin{align*}
		B &\leq \mE\bigg[\Bigl(\sum_{n = 1}^M \int_{t_{n-1}}^{t_n}\|\vG(\vu(s)) - \G(s;\vu(t_{n-1}))\|^2_{\vL^2}\|\ve^{n-1}\|^2_{\vL^2}\, ds\Bigr)^{1/2}\bigg] \\\nonumber
		&+ \mE\bigg[\Bigl(\sum_{n = 1}^M \int_{t_{n-1}}^{t_n}\|\G(s;\vu(t_{n-1})) - \G(s;\vu^{n-1})\|^2_{\vL^2}\|\ve^{n-1}\|^2_{\vL^2}\, ds\Bigr)^{1/2}\bigg]\\\nonumber
		&\leq \frac18 \mE\bigl[\max_{1\leq m \leq M} \|\ve^n\|^2_{\vL^2}\bigr] + C k^{2(1-\eps)}.
	\end{align*}
	Substituting the estimates from $A$ and $B$ into \eqref{eq3.27} we obtain
	\begin{align}
		&\frac12\mE\Bigl[\max_{1\leq m \leq M}\|\ve^m\|^2_{\vL^2}\Bigr] + \frac12\mE\Bigl[\sum_{n = 1}^M \|\ve^n - \ve^{n-1}\|^2_{\vL^2}\Bigr] + \mE\Bigl[k\sum_{n = 1}^M \|\nab\ve^n\|^2_{\vL^2}\Bigr]\\\nonumber
		&\leq Ck^{2(1-\eps)} + \frac14\mE\bigl[\max_{1\leq n \leq M}\|\ve^n\|^2_{\vL^2}\bigr],
	\end{align}
	which implies the desired estimate \eqref{eq310}. The proof is complete.
\end{proof}

\bigskip

In the next theorem, we will establish an $O(k^{1-\eps})$ convergence order for the pressure approximation $\{p^n\}$ of the time discretization. As stated, the pressure solution $p$ is only a distribution, we will derive the error estimate for between $\{p^n\}$ and the ``simulate" pressure solution $\tilde{p}$. Moreover, let $\delta = k$ in \eqref{simu_press}, then for $1 \leq n \leq M$ 
\begin{align*}
	\tilde{p}(t_n) = \frac{P(t_n) - P(t_{n-1})}{k}.
\end{align*}

\begin{theorem}\label{thm_semi_pressure}
	Under the assumptions of Theorem \ref{theorem_semi}. For $\eps>0$, there exists a constant $C= C(\vu_0,D_T)>0$, such that 
	\begin{align}\label{eq3.20}
		\mE\Bigl[k\sum_{n = 1}^M\bigl\|\tilde{p}(t_n) - p^n\bigr\|_{L^2}\Bigr]\leq C k^{1-\eps}.
	\end{align}
\end{theorem}
\smallskip
\begin{proof} 
	The proof is based on the well-known inf-sup (LBB) condition which is quoted below. There exists a constant $\beta >0$ such that
	\begin{align}\label{inf-sup}
		\beta \|\xi\|_{L^2} \leq \sup_{\pphi \in \vH^1_{per}(D)}\frac{\bigl(\xi,\div \pphi\bigr)}{\|\nab\pphi\|_{\vL^2}}\,\,\qquad\qquad\forall \xi \in L_{per}^2(D).
	\end{align}
	
	Subtracting \eqref{eq2.10a} to \eqref{milsteinscheme1} and using the notation $\ve^n:= \vu(t_n) - \vu^n$, we obtain
	\begin{align*}
		k\bigl(\div \pphi, \tilde{p}(t_n) - p^n\bigr) &=	\bigl(\ve^n - \ve^{n-1},\pphi\bigr) + k\bigl(\nab\ve^n,\nab\pphi\bigr) \\\nonumber
		&\qquad+ \int_{t_{n-1}}^{t_n} \bigl(\nab\bigl(\vu(s) - \vu(t_n)\bigr),\nab\pphi\bigr)\, ds \\\nonumber
		&\qquad- \Bigl(\int_{t_{n-1}}^{t_n}\bigl(\vG(\vu(s)) - \G(s;\vu^{n-1})\bigr)\, dW(s),\pphi\Bigr).
	\end{align*}
	
	By using \eqref{inf-sup}, for all $0\neq \pphi \in \vH_{per}^1(D)$, we have
	\begin{align}\label{eq3.21}
		\nonumber	k\beta\bigl\|\tilde{p}(t_n) - p^n\bigr\|_{L^2} &\leq \sup_{\pphi \in \vH^1_{per}(D)} \frac{k\Bigl(\div \pphi, \tilde{p}(t_n) - p^n\Bigr)}{\|\nab\pphi\|_{\vL^2}}\\\nonumber
		&\leq \frac{1}{\|\nab\pphi\|_{\vL^2}}\Bigl\{\bigl(\ve^n - \ve^{n-1},\pphi\bigr) + k\bigl(\nab\ve^n,\nab\pphi\bigr) \\
		&\quad+ \int_{t_{n-1}}^{t_n} \bigl(\nab\bigl(\vu(s) - \vu(t_n)\bigr),\nab\pphi\bigr)\, ds \\\nonumber
		&\quad- \Bigl(\int_{t_{n-1}}^{t_n}\bigl(\vG(\vu(s)) - \G(s;\vu^{n-1})\bigr)\, dW(s),\pphi\Bigr)\Bigr\}
	\end{align}
	
	Next, applying the summation $\sum_{n = 1}^M$ to \eqref{eq3.21} yields to
	\begin{align}\label{eq3.22}
		k\sum_{n = 1}^M \bigl\|\tilde{p}(t_n) - p^n\bigr\|_{L^2} &\leq \frac{1}{\beta\|\nab\pphi\|_{\vL^2}}\Bigl\{\bigl(\ve^M,\pphi\bigr) + \Bigl(k\sum_{n = 1}^M\nab\ve^n,\nab\pphi\Bigr) \\\nonumber
		&+ \sum_{n = 1}^M\int_{t_{n-1}}^{t_n} \bigl(\nab\bigl(\vu(s) - \vu(t_n)\bigr),\nab\pphi\bigr)\, ds \\\nonumber
		&-\Bigl(\sum_{n = 1}^M\int_{t_{n-1}}^{t_n} \bigl(\vG(\vu(s)) - \G(s;\vu^n)\bigr)\, dW(s),\pphi\Bigr)\Bigr\}\\\nonumber
		&:= {\tt I + II + III + IV}.
	\end{align}
	Now, we can estimate {\tt I, $\cdots$, IV} as follows.
	
	By using Schwarz inequality, Poincar\'e inequality and Theorem \ref{theorem_semi} we obtain
	\begin{align}
		\mE\bigl[{\tt I + II}\bigr] \leq C\mE\bigl[\|\nab\ve^M\|_{\vL^2} + k\sum_{n = 1}^M\|\nab\ve^n\|_{\vL^2}\bigr] \leq C\, k^{1 - \eps}.
	\end{align}
	
	Next, by using Lemma \ref{lemma3.1}$(ii)$ and martingale property of the It\^o integral we obtain
	\begin{align*}
		\mE\bigl[{\tt III}\bigr] &=  -\mE\Bigl[\frac{1}{\|\nab \pphi\|_{\vL^2}}\sum_{n = 1}^M\int_{t_{n-1}}^{t_n} \Bigl(\nab\Bigl(\vu(t_n) - \vu(s) -\int_{s}^{t_n} \vG(\vu(\xi))\, dW(\xi) \Bigr),\nab\pphi\Bigr)\, ds\Bigr]\\\nonumber
		&\qquad-\mE\Bigl[\frac{1}{\|\nab\pphi\|_{\vL^2}}\sum_{n = 1}^M\int_{t_{n-1}}^{t_n}\Bigl(\int_{s}^{t_n}\nab\vG(\vu(\xi))\, dW(\xi),\nab\pphi\Bigr)\, ds\Bigr]\\\nonumber
		&\leq \mE\Bigl[\sum_{n = 1}^M \int_{t_{n-1}}^{t_n} \Bigl\|\vu(t_n) - \vu(s) - \int_{s}^{t_n} \vG(\vu(\xi))\, dW(\xi)\Bigr\|_{\vH^1} \, ds\Bigr]\\\nonumber
		&\leq C\, k^{1-\eps}.
	\end{align*}
	
	In addition, due to the martingale property of the It\^o integral, $\mE\bigl[{\tt IV}\bigr] = 0$.
	
	Finally, substituting the estimates from ${\tt I, \cdots, IV}$ into \eqref{eq3.22} we obtain the desired estimate. The proof is complete.
\end{proof}

\begin{remark} Theorem \ref{thm_semi_pressure} provides an error estimate for the approximate pressure solution $\{p^n\}$ in the discrete $L^1(\Ome; \ell^1(0,T;\vL^2(D)))$-norm.  The error estimate in this norm is stronger than in the time-averaged norm that has been used in \cite{FPL2021, FQ,FL2020,FVO2022, LV2021,vo2022} to analyze the error estimates of the pressure approximation of the Euler-Maruyama method for the same equation. Indeed, we have
	\begin{align}
		\mE\Bigl[\Bigl\|P(t_m) - k\sum_{n = 1}^m p^n\Bigr\|_{L^2}\Bigr] &= \mE\Bigl[\Bigl\|k\sum_{n = 1}^m\Bigl(\frac{P(t_n) - P(t_{n-1})}{k} - p^n\Bigr) \Bigr\|_{L^2}\Bigr] \\\nonumber
		&\leq \mE\Bigl[k\sum_{n = 1}^m\bigl\|\tilde{p}(t_n) - p^n\bigr\|_{L^2}\Bigr] \leq C\, k^{1-\eps}.
	\end{align}
	In addition, the new interpretation $\tilde{p}$ gives more advantage in computation that will be addressed in Section \ref{sec_compute}.
\end{remark}	

\bigskip
\section{Fully discrete mixed finite element discretization}\label{section_fullydiscrete} In this section, we formulate and 
analyze the spatial approximations of Algorithm 1 by using the mixed finite element method.

\subsection{Standard mixed finite element method}
Let $\mathcal{T}_h$ be a quasi-uniform mesh of the domain $D \subset \mathbb{R}^2$ with mesh size $h > 0$. We consider the following mini finite element pair:
\begin{align*}
	\mH_h &= \bigl\{\vv_h \in {\bf C}(\overline{D}) \cap \vH^1_{per}(D);\,\vv_h \in [\mathcal{P}_{1b}(K)]^2\qquad\forall K \in \mathcal{T}_h\bigr\},\\
	L_h &= \bigl\{\psi_h \in C(\overline{D})\cap \in L^2_{per}(D)/\mathbb{R};\, \psi_h \in \mathcal{P}_1(K)\qquad\forall K \in \mathcal{T}_h\bigr\},
\end{align*}  
where $\mathcal{P}_{1}$ is the space of linear piece-wise polynomials on $K$. Moreover, $\mathcal{P}_{1b}$ denotes the space of buble linear polynomials on $K$. It is well-known that the mini finite element space pair $\mH_h$ and $L_h$ satisfies the Ladyzhenskaja-Babuska-Brezzi (LBB)  (or inf-sup condition) which is now quoted:  there exists $\beta_1>0$ such that
\begin{align}\label{inf-sup_discrete}
	\sup_{\pphi_h \in \mH_h} \frac{\bigl(\div \pphi_h,\psi_h\bigr)}{\|\nab\pphi_h\|_{\vL^2}} \geq \beta_1\|\psi_h\|_{L^2}\qquad\forall \psi_h\in L_h,
\end{align}
where the constant $\beta_1$ is independent of $h$ (and $k$). 

\smallskip
\noindent
\textbf{Algorithm 2} 

Let $\vu_h^0$ be a given $\mH_h$-valued random variable. Find $\displaystyle \bigl(\vu_h^{n+1},p_h^{n+1}\bigr) \in \mH_h\times L_h$ such that $\mP$-a.s.
\begin{align}
	\label{eq4.2}	\bigl(\vu_h^{n+1} &- \vu_h^n,\pphi_h\bigr) + \nu k \bigl(\nab \vu_h^{n+1},\nab\pphi_h\bigr) - k\bigl(p_h^{n+1},\div \pphi_h\bigr) \\\nonumber
	&\qquad\qquad= \bigl(\vG(\vu^n_h)\Delta W_{n} + \frac12 D\vG(\vu_h^n) \vG(\vu_h^n)[(\Delta W_{n})^2 - k],\pphi_h\bigr),\\
	&\bigl(\div\vu_h^n,\psi_h\bigr) =0,
\end{align}
for all $\pphi_h \in \mH_h$ and $\psi_h\in L_h$ and $\Delta W_n = W(t_{n+1}) - W(t_n)\sim \mathcal{N}(0,k)$.

\smallskip

Upon using the functional $\mathcal{G}$ defined in \eqref{eq3.2} we can rewrite \eqref{eq4.2} as follow
\begin{align}\label{eq4.4}
	\bigl(\vu_h^{n+1} - \vu_h^n,\pphi_h\bigr) + \nu k \bigl(\nab \vu_h^{n+1},\nab\pphi_h\bigr) &- k\bigl(p_h^{n+1},\div \pphi_h\bigr) \\\nonumber
	&= \Bigl(\int_{t_n}^{t_{n+1}}\G(s;\vu^n_h)\, dW(s), \pphi_h\Bigr).
\end{align}

Next, we define the following space
\begin{align*}
	\mV_h = \Bigl\{\pphi_h \in \mH_h;\, \bigl(\div \pphi_h, q_h\bigr) =0\qquad\forall q_h \in L_h\Bigr\}.
\end{align*}

We observe that $\mV_h \subset \mH_h$ and in general, $\mV_h$ is not a subspace of $\mV$.

Denote $\vQ_h: \vL^2_{per} \rightarrow \mV_h$ as the $L^2$-orthogonal projection, which satisfies 
\begin{align}\label{L2project}
	\bigl(\vv - \vQ_h\vv, \pphi_h\bigr) = 0\qquad\forall \pphi_h\in \mV_h.
\end{align}

In addition, we recall the following well-known interpolation estimates:
\begin{align}\label{L2inequ}
	\|\vv - \vQ_h\vv\|_{\vL^2} + h\|\nab(\vv - \vQ_h\vv)\|_{\vL^2} &\leq C h^2\|\vA\vv\|_{\vL^2}\qquad\forall \vv\in \mV\cap\vH^2(D),\\
	\|\vv - \vQ_h\vv\|_{\vL^2} &\leq Ch\|\nab\vv\|_{\vL^2}\qquad\forall \vv\in \mV\cap\vH^1(D).
\end{align} 

We also let $P_h: L^2_{per} \rightarrow L_h$ denote the $L^2$-orthogonal projection defined by   
\begin{align}
	\bigl(\psi - P_h\psi, q_h\bigr) = 0\qquad\forall q_h \in L_h.
\end{align}
It is well-known that there holds 
\begin{align}
	\|\psi - P_h\psi\|_{L^2} \leq Ch\|\nab \psi\|_{\vL^2}\qquad\forall \psi \in L^2_{per}(D)\cap H^1(D).
\end{align}

We state the following stability estimate for $\{\vu_h^n\}$. Its proof is similar to the proof of \cite[Lemma 3.1]{BCP12}, so we omit the proof to save space.

\begin{lemma}\label{stability_FEM}
	Let $\vu_h^0 \in L^{2}(\Ome;\mH_h)$ satisfying  $ \mE\bigl[\|\vu_h^0\|^{2}_{\vL^2}\bigr] \leq C$. Then, there exists a pair $\bigl\{\vu_h^{n},p^{n}_h\bigr\}_{n=1}^M \subset L^{2}(\Ome; \mH_h\times L_h)$ that solves Algorithm 2 and satisfies
	\begin{align}
		\mE\biggl[\max_{1\leq n \leq M}\|\vu_h^n\|^{2}_{\vL^2} + \nu k\sum_{n=1}^M \|\nab\vu_h^n\|^2_{\vL^2}\biggr] \leq C.
	\end{align}
	where $C = C(D_T,\vu_h^0)>0$.
\end{lemma}


\subsection{Error estimates for Algorithm 2}

In this part, we state error estimates for the finite element approximation velocity $\vu_h^n$ and the semi-discretization $\vu^n$ in Theorem \ref{theorem_fully} below. The proof of this theorem is similar to the proof of \cite[Theorem 4.2]{FPL2021} with a few modifications. However, for the sake of completeness, we also present that proof below. 

\begin{theorem}\label{theorem_fully} 
	Let $\vu_h^0 = \vQ_h \vu_0$.  Assume that $\vG$ satisfies ${\bf (A1), (A2), (A3)}$ and $\vu_0 \in L^{2}(\Ome;\mV)$. Let $\{(\vu^n, p^n)\}$ and $\{\vu^n_h,p^n_h\}$ be the velocity and pressure approximations generated by Algorithm 1 and Algorithm 2, respectively. Then there holds
	\begin{align}\label{fully_error}
		\bigl(\mE\bigl[\max_{1\leq n \leq M}|\vu^n - \vu_h^n\|^{2}_{\vL^2}\bigr]\bigr)^{\frac12} &+ \Bigl(\mE\Bigl[ k \sum_{n=1}^M\|\nab(\vu^n - \vu_h^n)\|^2_{\vL^2}\Bigr]\Bigr)^{\frac12} \leq C\,h,
	\end{align}
	where $C= C(\vu_0,D_T)>0$ is independent of $k$ and $h$.
\end{theorem}	

\begin{proof}
	First, define $\ve_{\vu}^n := \vu^n - \vu^n_h$ and $e^n_p := p^n - p^n_h$. By subtracting \eqref{milsteinscheme3} to \eqref{eq4.4} we obtain the following error equations
	\begin{align}\label{eq4.11}
		\bigl(\ve_{\vu}^{n+1} - \ve_{\vu}^n,\pphi_h\bigr) +  k\bigl(\nab\ve_{\vu}^{n+1},&\nab\pphi_h\bigr) - k\bigl(e^n_p,\div \pphi_h\bigr) \\\nonumber
		&= \Bigl(\int_{t_{n}}^{t_{n+1}} \bigl(\G(s;\vu^n) - \G(s;\vu^n_h)\bigr)\, dW(s),\pphi_h\Bigr),\\
		\bigl(\div \ve_{\vu}^{n+1}, q_h\bigr) &=0.
	\end{align}
	Choosing $\pphi_h = \vQ_h\ve_{\vu}^{n+1}\in \mV_h$ in \eqref{eq4.11} and using the orthogonality of $\vQ_h$, we obtain
	\begin{align}\label{eq4.13}
		\bigl(\vQ_h\ve_{\vu}^{n+1} &- \vQ_h\ve_{\vu}^n,\vQ_h\ve_{\vu}^{n+1}\bigr) + k \|\nab\vQ_h\ve_{\vu}^{n+1}\|^2_{\vL^2} \\\nonumber
		&=- k\bigl(\nab(\ve_{\vu}^{n+1}- \vQ_h\ve_{\vu}^{n+1}),\nab\vQ_h\ve_{\vu}^{n+1}\bigr) + k\bigl(e^n_p,\div \vQ_h\ve_{\vu}^{n+1}\bigr) \\\nonumber
		&\qquad+ \Bigl(\int_{t_{n}}^{t_{n+1}} \bigl(\G(s;\vu^n) - \G(s;\vu^n_h)\bigr)\, dW(s),\vQ_h\ve_{\vu}^{n+1} - \vQ_h\ve_{\vu}^n\Bigr)\\\nonumber
		&\qquad +  \Bigl(\int_{t_{n}}^{t_{n+1}} \G(s;\vu^n) - \G(s;\vu^n_h)\bigr)\, dW(s), \vQ_h\ve_{\vu}^n\Bigr)\\\nonumber
		&:= {\tt I + II + III + IV}.
	\end{align}
	
	We now can estimate the right-hand of \eqref{eq4.13} as follows.
	
	First, by using \eqref{L2project} and then \eqref{L2inequ} we have
	\begin{align*}
		{\tt I} &= -k\bigl(\nab(\vu^{n+1}- \vQ_h\vu^{n+1}),\nab\vQ_h\ve_{\vu}^{n+1}\bigr) \\\nonumber
		&\leq Ck\|\nab(\vu^{n+1} - \vQ_h\vu^{n+1})\|^2_{\vL^2} + \frac{ k}{4} \|\nab\vQ_h\ve_{\vu}^{n+1}\|^2_{\vL^2}\\\nonumber
		&\leq Ckh^2\|\vu^{n+1}\|^2_{\vH^2} + \frac{k}{4} \|\nab\vQ_h\ve_{\vu}^{n+1}\|^2_{\vL^2}.
	\end{align*}
	
	Next, by using the fact that $\bigl(q_h,\div\vQ_h\ve_{\vu}^{n+1}\bigr) = 0$ for all $q_h \in L_h$ we have
	\begin{align*}
		{\tt II} = k\bigl(p^{n+1},\div \vQ_h\ve_{\vu}^{n+1}\bigr) &= k\bigl(p^{n+1} - P_h p^{n+1},\div\vQ_h\ve_{\vu}^{n+1}\bigr)\\\nonumber
		&\leq k\|p^{n+1} - P_h p^{n+1}\|^2_{L^2} + \frac{k}{4}\|\nab\vQ_h\ve_{\vu}^{n+1}\|^2_{\vL^2}\\\nonumber
		&\leq Ckh^2\|\nab p^{n+1}\|^2_{\vL^2} + \frac{k}{4}\|\nab\vQ_h\ve_{\vu}^{n+1}\|^2_{\vL^2}.
	\end{align*}
	
	Now, we observe that $\mE[{\tt IV}] = 0$ due to the martingale property of the It\^o integral. It is left to estimate ${\tt III}$ as below. By using the It\^o isometry and Lemma \ref{lemma3.1}$(i)$ we obtain
	\begin{align*}
		\mE\bigl[{\tt III}\bigr] &\leq \mE\Bigl[\Bigl\|\int_{t_{n}}^{t_{n+1}} \bigl(\G(s;\vu^n) - \G(s;\vu_h^n)\bigr)\, dW(s)\Bigr\|^2_{\vL^2}\Bigr] + \frac14 \mE\bigl[\|\vQ_h(\ve_{\vu}^{n+1} - \ve_{\vu}^n) \|^2_{\vL^2}\bigr]\\\nonumber
		&=  \mE\Bigl[\int_{t_{n}}^{t_{n+1}} \big\|\G(s;\vu^n) - \G(s;\vu_h^n)\bigr\|^2_{\vL^2}\, ds\Bigr] + \frac14 \mE\bigl[\|\vQ_h(\ve_{\vu}^{n+1} - \ve_{\vu}^n) \|^2_{\vL^2}\bigr]\\\nonumber
		&\leq Ck\mE\bigl[\|\vQ_h\ve_{\vu}^n\|^2_{\vL^2}\bigr] + Ckh^4\mE\bigl[\|\vu^{n}\|^2_{\vH^2}\bigr]+\frac14 \mE\bigl[\|\vQ_h(\ve_{\vu}^{n+1} - \ve_{\vu}^n) \|^2_{\vL^2}\bigr].
	\end{align*}
	
	The left-hand side of \eqref{eq4.13} can be analyzed by using the identity $2a(a-b) = a^2 - b^2 + (a-b)^2$ as follows.
	\begin{align*}
		\bigl(\vQ_h\ve_{\vu}^{n+1} &- \vQ_h\ve_{\vu}^n,\vQ_h\ve_{\vu}^{n+1}\bigr)  = \frac{1}{2}\bigl[\|\vQ_h\ve_{\vu}^{n+1}\|^2_{\vL^2} - \|\vQ_h\ve_{\vu}^n\|^2_{\vL^2}\bigr] + \frac12 \|\vQ_h\ve_{\vu}^{n+1} - \vQ_h\ve_{\vu}^n\|^2_{\vL^2}.
	\end{align*}
	Therefore, substituting the estimates from {\tt I, $\cdots$, IV} into \eqref{eq4.13} and then absorbing the like-terms from the right side to the left side we obtain
	\begin{align}\label{eq4.14}
		&\frac12\mE\bigl[\|\vQ_h\ve_{\vu}^{n+1}\|^2_{\vL^2} - \|\vQ_h\ve_{\vu}^n\|^2_{\vL^2} +\frac14\|\vQ_h(\ve_{\vu}^{n+1} - \ve_{\vu}^n)\|^2_{\vL^2}\bigr] + \frac{k}{2}\mE\bigl[\|\nab\vQ_h\ve_{\vu}^{n+1}\|^2_{\vL^2}\bigr]\\\nonumber
		&\leq Ckh^2\mE\bigl[\|\vu^{n+1}\|^2_{\vH^2}\bigr] + Ckh^2\mE\bigl[\|\nab p^{n+1}\|^2_{\vL^2}\bigr] + Ck\mE\bigl[\|\vQ_h\ve_{\vu}^n\|^2_{\vL^2}\bigr].
	\end{align}
	
	Next, applying the summation $\sum_{n=0}^m$ to \eqref{eq4.14} for any $0\leq m \leq M-1$ we obtain
	\begin{align}\label{eq4.15}
		&\frac12\mE\bigl[\|\vQ_h\ve_{\vu}^{m+1}\|^2_{\vL^2}\bigr] +\frac14\sum_{n=0}^m\mE\bigl[\|\vQ_h(\ve_{\vu}^{n+1} - \ve_{\vu}^n)\|^2_{\vL^2}\bigr] + \frac{k}{2}\sum_{n=0}^m\mE\bigl[\|\nab\vQ_h\ve_{\vu}^{n+1}\|^2_{\vL^2}\bigr]\\\nonumber
		&\leq Ch^2k\sum_{n=0}^m\mE\bigl[\|\vu^{n+1}\|^2_{\vH^2}\bigr] + Ch^2k\sum_{n=0}^M\mE\bigl[\|\nab p^{n+1}\|^2_{\vL^2}\bigr] + Ck\sum_{n=0}^m\mE\bigl[\|\vQ_h\ve_{\vu}^n\|^2_{\vL^2}\bigr].
	\end{align}
	By using the discrete Gronwall inequality to \eqref{eq4.15} and then using the stability estimates in Lemma \ref{stability_mean}, we obtain
	\begin{align}\label{eq4.16}
		&\frac12\mE\bigl[\|\vQ_h\ve_{\vu}^{m+1}\|^2_{\vL^2}\bigr] +\frac14\sum_{n=0}^m\mE\bigl[\|\vQ_h(\ve_{\vu}^{n+1} - \ve_{\vu}^n)\|^2_{\vL^2}\bigr] + \frac{k}{2}\sum_{n=0}^m\mE\bigl[\|\nab\vQ_h\ve_{\vu}^{n+1}\|^2_{\vL^2}\bigr]\\\nonumber
		&\leq Ch^2 \exp(CT).
	\end{align}
	By applying $\max_{0\leq m \leq M-1}$ to \eqref{eq4.16}, we obtain
	\begin{align}\label{eq4.17}
		\max_{1\leq m \leq M}\mE\bigl[\|\vQ_h\ve_{\vu}^m\|^2_{\vL^2}\bigr] + \mE\Bigl[ k\sum_{n=1}^M\|\nab\vQ_h\ve_{\vu}^{n}\|^2_{\vL^2}\Bigr] \leq Ch^2.
	\end{align}
	Next, we use the estimate \eqref{eq4.17} to prove the desired estimate \eqref{fully_error}. To do that, we follow the techniques that have been used in the proof of Theorem \ref{theorem_semi} as follows. 
	
	First, applying $\sum_{n=1}^m$ (for $1\leq m \leq M$), $\max_{1\leq m \leq M}$, $\mE[\cdot]$ to \eqref{eq4.13}, respectively, we obtain
	\begin{align}
		&\frac12\mE\Bigl[\max_{1\leq m \leq M}\|\vQ_h\ve_{\vu}^m\|^2_{\vL^2} + \sum_{n=1}^M\|\vQ_h\ve_{\vu}^n - \vQ_h\ve_{\vu}^{n-1}\|^2_{\vL^2} + 2 k \sum_{n=1}^M \|\nab\vQ_h\ve_{\vu}^n\|^2_{\vL^2}\Bigr]\\\nonumber
		&\leq Ch^2\mE\Bigl[k\sum_{n=1}^M\|\vu^n\|^2_{\vH^2} + k\sum_{n = 1}^M\|\nab p^n\|^2_{\vL^2}\Bigr] + \frac14\mE\Bigl[k\sum_{n = 1}^M\|\nab\vQ_h\ve_{\vu}^n\|^2_{\vL^2}\Bigr]\\\nonumber
		&\leq \mE\Bigl[\sum_{n = 1}^{M}\Bigl\|\int_{t_{n-1}}^{t_n}\bigl(\G(s;\vu^{n-1}) - \G(s;\vu_h^{n-1})\bigr)\, dW(\xi)\Bigr\|^2_{\vL^2} \Bigr]\\\nonumber
		&\qquad+\mE\Bigl[\frac14\sum_{n = 1}^M \|\vQ_h\ve_{\vu}^{n} - \vQ_h\ve_{\vu}^{n-1}\|^2_{\vL^2}\Bigr]\\\nonumber
		&\qquad+ \mE\Bigl[\Bigl(\sum_{n = 1}^M \int_{t_{n-1}}^{t_n}\|\G(s;\vu^{n-1}) - \G(s;\vu_h^{n-1})\|^2_{\vL^2}\|\vQ_h\ve_{\vu}^{n-1}\|^2_{\vL^2}\, ds\Bigr)^{1/2}\Bigr]\\\nonumber
		&\leq Ch^2 + \frac14\mE\bigl[\max_{1\leq m \leq M}\|\vQ_h\ve_{\vu}^m\|^2_{\vL^2}\bigr],
	\end{align}
	which implies that
	\begin{align}
		\mE\bigl[\max_{1\leq n \leq M}\|\vQ_h\ve_{\vu}^n\|^2_{\vL^2}\bigr] + \mE\Bigl[k\sum_{n = 1}^M\|\nab\vQ_h\ve_{\vu}^n\|^2_{\vL^2}\Bigr] \leq Ch^2.
	\end{align}
	Finally, the proof is completed by a simple triangular inequality process.
\end{proof}

\smallskip

Next, we also present the error estimate of the pressure approximations $\{p_h^n\}$ and  $\{p^n\}$ from Algorithms 1 and 2, respectively. The proof of this theorem is similar to the one of Theorem \ref{thm_semi_pressure} so we leave it as an exercise for the reader to try.
\begin{theorem}\label{theorem_fully_pressure}
	Under the same assumptions of Theorem \ref{theorem_fully}. There exists a constant $C = C(\beta_1,\vu_0,D_T)>0$ such that
	\begin{align*}
		\mE\Bigl[k\sum_{n = 1}^M\|p^n - p^n_h\|_{L^2}\Bigr] \leq C\, h.
	\end{align*}
\end{theorem}

\medskip

Finally, we are ready to state the global error estimates of the fully discrete velocity and pressure approximations. 
\begin{theorem}\label{thm_global}
	Let $\eps > 0$. Under the assumptions of Theorem \ref{theorem_semi} and Theorem \ref{theorem_fully}, there exist a constant $C = C(D_T, \vu_0, \beta, \beta_1)>0$ such that
	\begin{align*}
		\bigl(\mE\bigl[\max_{1\leq n \leq M}|\vu(t_n) - \vu_h^n\|^{2}_{\vL^2}\bigr]\bigr)^{\frac12} + \Bigl(\mE\Bigl[\nu k \sum_{n=1}^M\|\nab(\vu(t_n) - \vu_h^n)\|^2_{\vL^2}\Bigr]\Bigr)^{\frac12} 
		&\leq C\,\bigl(k^{1-\eps} + h\,\bigr),\\\nonumber
		\mE\Bigl[k\sum_{n = 1}^M\bigl\|\tilde{p}(t_n) - p^n_h\bigr\|_{L^2}\Bigr] &\leq C\, \bigl(k^{1-\eps} + h\bigr).
	\end{align*}
\end{theorem}

\section{Computational experiments}\label{sec_compute}

In this section, we present numerical tests to verify our convergence order $O(k^{1-\eps} + h)$ in Theorem \ref{thm_global}. In all our experiments   we set $D = (0,1)^2\subset \mathbb{R}^2$, $T =1, \nu =1$, the body force is $\vf = (f_1,f_2)$ with
\begin{align*}
	f_1(x,y) &= \pi\cos(t)\sin(2\pi y)\sin(\pi x)\sin(\pi x)
	- 2\pi^3\sin(t)\sin(2\pi y)(2\cos(2\pi x)-1) \\\nonumber
	&\qquad- \pi\sin(t)\sin(\pi x)\sin(\pi y),\\\nonumber
	f_2(x,y) &= -\pi\cos(t)\sin(2\pi x)\sin(\pi y)\sin(\pi y)
	- 2\pi^3\sin(t)\sin(2\pi x)(1-2\cos(2\pi y)) \\\nonumber
	&\qquad+ \pi\sin(t)\cos(\pi x)\cos(\pi y).
\end{align*}
We choose $W(t)$ in \eqref{eq1.1} to be a $\mathbb{R}$-valued Wiener process that is simulated by the minimal time step size $k_0 = 1/2048$ and the number of samples $J = 300$. We use the standard Monte Carlo method to compute the expectation. We take $\vG(\vu) = \alpha\vu$ for all $\vu \in \mV$ and $\alpha>0$, which gives $D\vG(\vu)(\vv) = \alpha\vv$ $\forall \vv\in \mV$. In addition, we use the Mini finite element method for the spatial discretization and the homogeneous Dirichlet boundary condition is imposed on $\vu$. 

We implement Algorithm 2 and compute the errors of the velocity and pressure approximations in the specified norms below. Since the exact solutions are unknown, the errors are computed between the computed solution $(\vu_h^n(\omega_j), p_h^n(\omega_j))$ and a reference solution $(\vu^n_{ref}(\omega_j), {{p}}^n_{ref}(\omega_j))$ (specified later) at the $\omega_j$-th sample.

Furthermore, to evaluate errors in strong norms, we use the following numerical 
integration formulas:
\begin{align*}
	L^2_{\omega}L^{\infty}_tL^2_x(\vu)&:=	\bigl(\mE\bigl[\max_{1\leq n \leq M}\|\vu(t_n) - \vu^n_h\|^2_{\vL^2}\bigr]\bigr)^{1/2} \\
	&\approx \Bigl(\frac{1}{J}\sum_{j= 1}^J\bigl(\max_{1\leq n \leq M} \|\vu^n_{ref}(\omega_j) - \vu^n_h(\omega_j)\|^2_{\vL^2}\bigr)\Bigr)^{1/2},\\
	L^2_{\omega}L^2_tH^1_x(\vu)&:= \Bigl(\mE\Bigl[k\sum_{n = 1}^M \|\vu(t_n) - \vu^n_h\|^2_{\vH^1}\Bigr]\Bigr)^{1/2} \\
	&\approx \Bigl(\frac{1}{J}\sum_{j= 1}^J\Bigl(k\sum_{n = 1}^M\|\vu^n_{ref}(\omega_j) - \vu_h^n(\omega_j)\|^2_{\vH^1}\Bigr)\Bigr)^{1/2},\\
	L^1_{\omega}L^1_tL^2_x(p) &: = \mE\Bigl[k\sum_{n = 1}^M\|\tilde{p}(t_n) - p^n_h\|_{L^2}\Bigr]\\
	&\approx \frac{1}{J}\sum_{j= 1}^J \bigg(k\sum_{n = 1}^M \|p^n_{ref}(\omega_j) - p^n_h(\omega_j)\|_{L^2}\bigg).
\end{align*}

{\bf Test 1.} In the first test, we want to verify the convergence order of the time discretization in Algorithm 1. To do that, we run Algorithm 2 to compute $\{(u^n_h,p^n_h)\}$ with a fixed mesh size $h = 1/40$ and vary the time step size by choosing $k = 2^{\ell} k_0$ for $\ell \in \mathbb{N}$ and the reference solutions $\{(\vu^n_{ref}, p^n_{ref})\}$ with $k_{ref} = k/2$ (i.e. we approximate the errors by comparing the numerical solutions in two consecutive time discretizations). The result errors are shown in Table \ref{table1}. The numerical results verify convergence order of almost $1$ for both velocity and pressure approximations as predicted by our error estimate results in Theorem \ref{theorem_semi}. 

\begin{table}[tbhp]
	\begin{center}
		\begin{tabular}{ |c|c|c|c|c|c|c|}
			\hline
			\bf $k$ & $L^2_{\omega}L^2_tH^1_x(\vu)$ & \mbox{order} & $L^2_{\omega}L^{\infty}_tL^2_x(\vu)$ &  \mbox{order}&$ L^1_{\omega}L^1_tL^2_x(p)$&order\\
			\hline 
			$1/64$  & 1.035310&  & 0.0316739&&2.00612&\\
			\hline
			$1/128$  & 0.611458 &0.7597  & 0.0207695&0.6088&1.13635&0.8200\\
			\hline
			$1/256$  & 0.341012 & 0.8430 & 0.0122318&0.7640&0.61963&0.8749\\
			\hline
			$1/512$  & 0.177547 & 0.9416 & 0.00663272&0.8830&0.321795&0.9445\\
			\hline
			$1/1024$& 0.0928572 &0.9352  &0.00359361 & 0.8842&0.172365 &0.9007\\
			\hline
		\end{tabular}
		\caption{The time discretization errors and convergence order of the computed velocity $\{ {\bf u}^n_h\}$ and pressure $\{p^n_h\}$ with $\alpha = \frac{1}{2}$.}
		\label{table1}
	\end{center}
\end{table}

\begin{figure}
\begin{center}
		\includegraphics[scale=0.15]{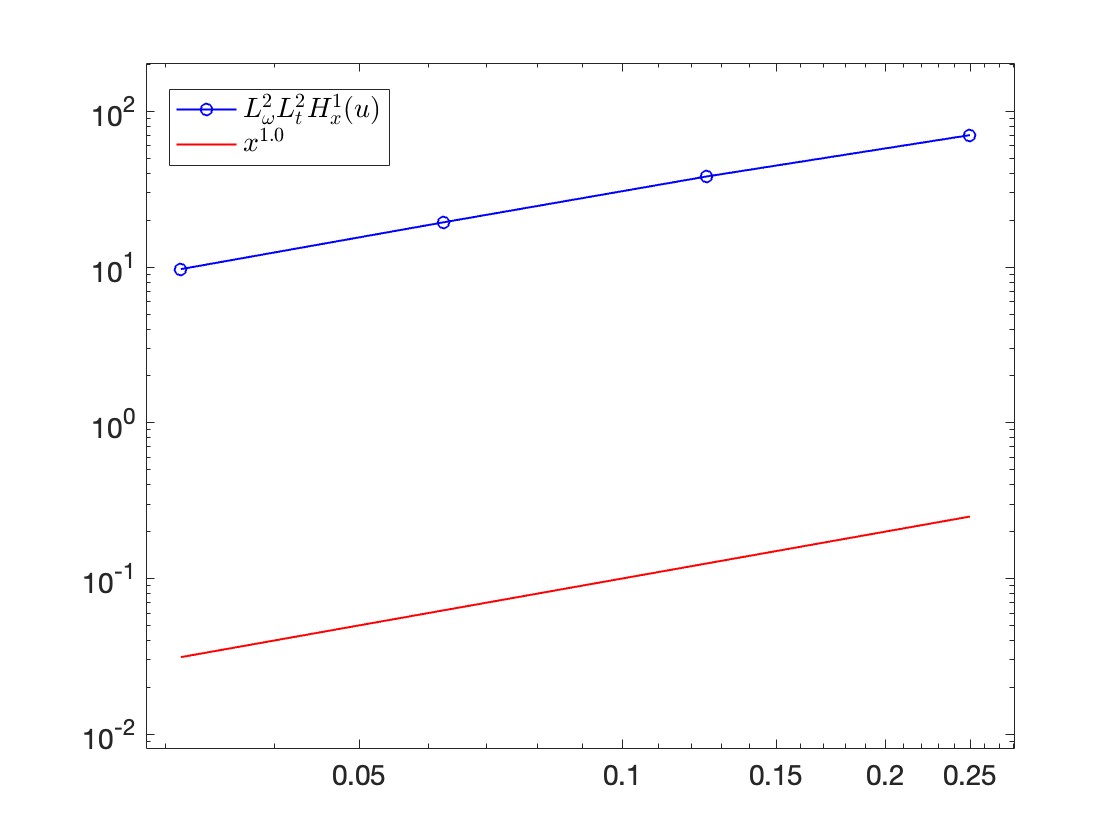}
	\includegraphics[scale=0.15]{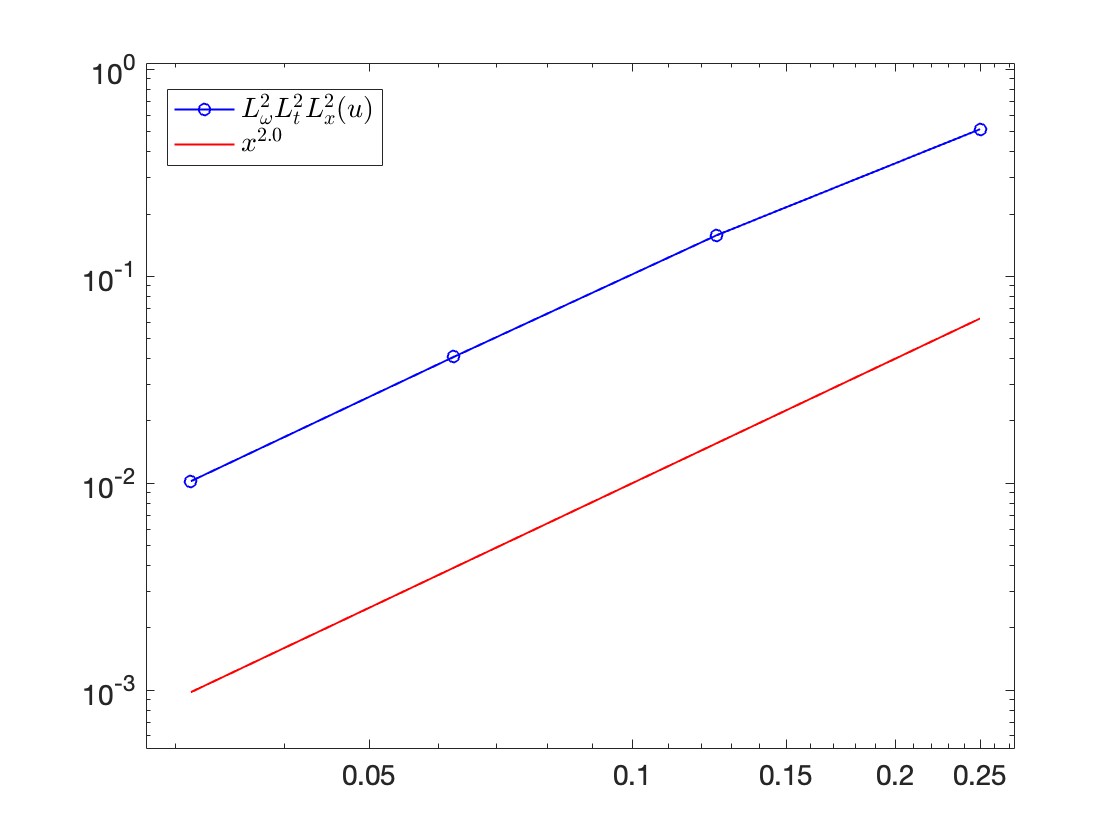}
	\caption{Plots of the spatial discretization errors and convergence order of the computed velocity $\{ {\bf u}^n_h\}$ with $\alpha = \frac{1}{2}$.}\label{fig5.1}
\end{center}
\end{figure}

\begin{figure}
	\begin{center}
		\includegraphics[scale=0.15]{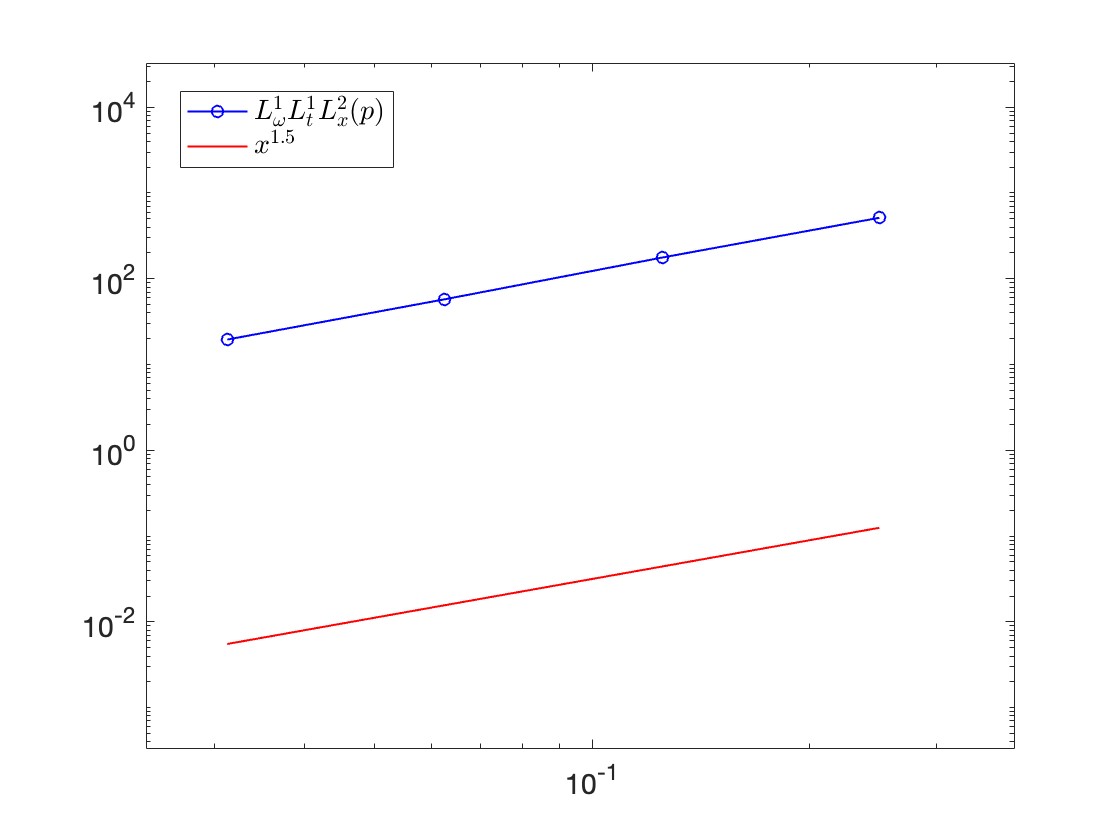}
		\caption{Plots of the spatial discretization errors and convergence order of the computed pressure $\{ {p}^n_h\}$ with $\alpha = \frac{1}{2}$.}\label{fig5.2}
	\end{center}
\end{figure}

{\bf Test 2.} In this test, we verify the spatial discretization convergence order of the mixed finite element method in Algorithm 2. To do that, we implement Algorithm 2 with a fixed time step size $k = 1/256$ and choose different mesh sizes $h$. The reference solutions $\{(\vu^n_{ref}, p^n_{ref})\}$ are computed by choosing $h_{ref} = h/2$. The result errors are displayed in Figures \ref{fig5.1} and \ref{fig5.2}. The numerical results in Figure \ref{fig5.1} (left) verify the first-order convergence of the velocity approximation in the $L^2_{\omega}L^2_tH^1_x$-norm as predicted in Theorem \ref{theorem_fully}. In addition, we also observe a second-order convergence (see Figure \ref{fig5.1} (right))for the computed velocity in the  $L^2_{\omega}L^2_tL^2_x$-norm which is better than our prediction in Theorem \ref{theorem_fully}. Furthermore, Figure \ref{fig5.2} also gives a convergence order of $1.5$ for the computed pressure that is slightly better than our predicted theory in Theorem \ref{theorem_fully_pressure}.



\section{Conclusion and discussion} We have established a convergence order of at most $1$ for the time-discretization of the stochastic Stokes equations. The paper proposes a new approach with analysis for the time-discretization based on the classical Milstein method, which is very popular in stochastic differential equations. However, the paper just considered the stochastic Stokes equations with a divergence-free multiplicative noise. This narrow assumption on the noise is due to the technical Lemma \ref{lemma3.1} which is unclear how to prove in the case of a non-divergence-free noise. In future work, we will extend the results of this paper to a more general multiplicative noise.

\medskip


\bibliographystyle{abbrv}
\bibliography{references}

\end{document}